\documentclass[twocolumn]{autart}

\usepackage{graphicx}
\usepackage{amssymb,amstext}
\usepackage{amsmath}
\usepackage{epsfig}
\usepackage{verbatim}
\usepackage[numbers,sort&compress]{natbib}

\newtheorem{theorem}{Theorem}[section]

\theoremstyle{definition}

\newtheorem{conjecture}{Conjecture}[section]

\theoremstyle{remark}
\newtheorem{remark}{Remark}

\newenvironment{proof}[1][Proof]{\begin{trivlist}
\item[\hskip \labelsep {\bfseries #1}]}{\end{trivlist}}

\renewcommand{\qed}{\nobreak \ifvmode \relax \else
      \ifdim\lastskip<1.5em \hskip-\lastskip
      \hskip1.5em plus0em minus0.5em \fi \nobreak
      \vrule height0.3em width0.5em depth0.25em\fi}

\begin{document}

\begin{frontmatter}
\title{Coherent-Classical Estimation for Linear Quantum Systems\thanksref{footnoteinfo}}

\thanks[footnoteinfo]{Preliminary versions of some of the results of this paper were presented at the $3^{rd}$ Australian Control Conference (AUCC), November 2013, Perth, Australia, and the $53^{rd}$ IEEE Conference on Decision and Control (CDC), December 2014, Los Angeles, CA, USA. Corresponding author: Shibdas Roy. Tel.~+44 (0) 24761 50635.}

\author{Shibdas Roy}$^{a}$\ead{roy\_shibdas@yahoo.co.in},
\author{Ian R. Petersen}$^{b}$\ead{i.r.petersen@gmail.com},
\author{Elanor H. Huntington}$^{b,c}$\ead{elanor.huntington@anu.edu.au}

\address{$^a$Department of Physics, University of Warwick, Coventry, United Kingdom}
\address{$^b$Research School of Engineering, Australian National University, Canberra, Australia}
\address{$^c$Australian Research Council Centre of Excellence for Quantum Computation and Communication Technology, Australia}

\begin{abstract}
We study a coherent-classical estimation scheme for a class of linear quantum systems, where the estimator is a mixed quantum-classical system that may or may not involve coherent feedback. We show that when the quantum plant or the quantum part of the estimator (coherent controller) is an annihilation operator only system, coherent-classical estimation without coherent feedback can provide no improvement over purely-classical estimation. Otherwise, coherent-classical estimation without feedback can be better than classical-only estimation for certain homodyne detector angles, although the former is inferior to the latter for the best choice of homodyne detector angle. Moreover, we show that coherent-classical estimation with coherent feedback is no better than classical-only estimation, when both the plant and the coherent controller are annihilation operator only systems. Otherwise, coherent-classical estimation with coherent feedback can be superior to purely-classical estimation, and in this case, the former is better than the latter for the optimal choice of homodyne detector angle.
\end{abstract}
\end{frontmatter}\vspace*{-4mm}

\vspace*{-7mm}
\section{Introduction}\vspace*{-2mm}
Estimation and control problems for quantum systems are of significant interest \cite{WM1,YK1,YK2,NY,JNP,NJP,GGY,MP1,MP2,GJN,IRP1}. An important class are linear quantum systems \cite{WM1,NY,JNP,NJP,GGY,GJN,GZ,WD,NJD,RP,IRP3}, that describe quantum optical devices such as optical cavities \cite{WM2,BR}, linear quantum amplifiers \cite{GZ}, and finite bandwidth squeezers \cite{GZ}. Coherent feedback control for linear quantum systems has been studied, where the feedback controller is also a quantum system \cite{JNP,NJP,MP1,MP2,WM3,SL}. A related coherent-classical estimation scheme was introduced by the authors in Refs. \cite{IRP2,RPH}, where the estimator has a classical part, which yields the desired final estimate, and a quantum part, which may involve coherent feedback. This is different from the quantum observer studied in Ref. \cite{MJ}. A quantum observer is a purely quantum system, that gives a quantum estimate of a variable for a quantum plant. By contrast, a coherent-classical estimator is a mixed quantum-classical system, that yields a classical estimate of a variable for a quantum plant.\vspace*{-3mm}

In this paper, we elaborate and build on the results of the conference papers Refs. \cite{RPH,IRP2} to present two key theorems, propose three relevant conjectures, and illustrate our findings with several examples. We show that a coherent-classical estimator without feedback, where either of the plant and the coherent controller is a physically realizable annihilation operator only system, it is not possible to get better estimates than the corresponding purely-classical estimator. Otherwise it is possible to get better estimates in certain cases. But we observe in examples that for the optimal choice of the homodyne angle, classical-only estimation is always superior. Moreover, we demonstrate that a coherent-classical estimator with coherent feedback can provide with higher estimation precision than classical-only estimation. This is possible only if either of the plant and the controller can not be defined purely using annihilation operators. Furthermore, if there is any improvement with the coherent-classical estimator (with feedback) over purely-classical estimation, we see in examples that the latter is always inferior for the optimal choice of the homodyne angle.\vspace*{-3mm}

The paper is structured as follows. Section \ref{sec:lqs_pr} introduces the class of linear quantum systems considered here and discusses physical realizability for such systems. Section \ref{sec:pce} formulates the problem of optimal purely-classical estimation. In Section \ref{sec:cce}, we formulate the optimal coherent-classical estimation problem without coherent feedback, and present our first theorem and two conjectures supported by examples. Section \ref{sec:ccef} discusses the coherent-classical estimation scheme involving coherent feedback and lays down our second theorem and another conjecture with pertinent examples. Finally, Section \ref{sec:conc} concludes the paper with relevant summarizing remarks.

\vspace*{-3mm}
\section{Linear Quantum Systems}\label{sec:lqs_pr}
The class of linear quantum systems we consider here are described by the quantum stochastic differential equations (QSDEs) \cite{JNP,GJN,IRP1,IRP2,SP}:

\vspace*{-8mm} \small
\begin{equation}\label{eq:lqs_1}
\begin{split}
\left[\begin{array}{c}
da(t)\\
da(t)^{\#}
\end{array}\right] &= F \left[\begin{array}{c}
a(t)\\
a(t)^{\#}
\end{array}\right] dt + G \left[\begin{array}{c}
d\mathcal{A}(t)\\
d\mathcal{A}(t)^{\#}
\end{array}\right];\\
\left[\begin{array}{c}
d\mathcal{A}^{out}(t)\\
d\mathcal{A}^{out}(t)^{\#}
\end{array}\right] &= H \left[\begin{array}{c}
a(t)\\
a(t)^{\#}
\end{array}\right] dt + K \left[\begin{array}{c}
d\mathcal{A}(t)\\
d\mathcal{A}(t)^{\#}
\end{array}\right],
\end{split}
\end{equation}
\normalsize \vspace*{-6mm}
where
\vspace*{-7mm} \small
\begin{equation}\label{eq:lqs_2}
\begin{split}
F &= \Delta(F_1,F_2), \qquad G = \Delta(G_1,G_2),\\
H &= \Delta(H_1,H_2), \qquad K = \Delta(K_1,K_2).
\end{split}
\end{equation}
\normalsize \vspace*{-6mm}

Here, $a(t) = [a_1(t) \hdots a_n(t)]^T$ is a vector of annihilation operators. The adjoint $a_i^{*}$ of the operator $a_i$ is called a creation operator. The notation $\Delta(F_1,F_2)$ denotes the matrix $\left[\begin{array}{cc} F_1 & F_2\\
F_2^{\#} & F_1^{\#}
\end{array}\right]$. Also, $F_1$, $F_2 \in \mathbb{C}^{n \times n}$, $G_1$, $G_2 \in \mathbb{C}^{n \times m}$, $H_1$, $H_2 \in \mathbb{C}^{m \times n}$, and $K_1$, $K_2 \in \mathbb{C}^{m \times m}$. Moreover, $^{\#}$ denotes the adjoint of a vector of operators or the complex conjugate of a complex matrix. Furthermore, $^\dagger$ denotes the adjoint transpose of a vector of operators or the complex conjugate transpose of a complex matrix. In addition, $\mathcal{A}(t) = [\mathcal{A}_1(t) \hdots \mathcal{A}_m(t)]^T$ is a vector of external independent quantum field operators and $\mathcal{A}^{out}(t)$ is the corresponding vector of output field operators.

\vspace*{-2mm}
\begin{theorem}\label{thm:phys_rlz}
(See \cite{SP,IRP3}) A complex linear quantum system of the form (\ref{eq:lqs_1}), (\ref{eq:lqs_2}) is physically realizable, if and only if there exists a complex commutation matrix $\Theta = \Theta^\dagger$ satisfying the following commutation relation
\vspace*{-8mm}\small
\begin{equation}\label{eq:comm_rel1}
\begin{split}
\Theta &=\left[\left[\begin{array}{c}
a\\
a^{\#}
\end{array}\right], \left[\begin{array}{c}
a\\
a^{\#}
\end{array}\right]^\dagger\right]\\
&= \left[\begin{array}{c}
a\\
a^{\#}
\end{array}\right]\left[\begin{array}{c}
a\\
a^{\#}
\end{array}\right]^\dagger - \left(\left[\begin{array}{c}
a\\
a^{\#}
\end{array}\right]^\# \left[\begin{array}{c}
a\\
a^{\#}
\end{array}\right]^T\right)^T,
\end{split}
\end{equation}
\normalsize \vspace*{-6mm}

such that
\vspace*{-1.2cm}\small
\begin{equation}
\begin{split}
F\Theta + \Theta F^\dagger + GJG^\dagger &= 0,\\
G &= -\Theta H^\dagger J,\\
K &= I,
\end{split}
\end{equation}
\normalsize \vspace*{-6mm}

where $J = \left[\begin{array}{cc}
I & 0\\
0 & -I
\end{array}\right]$.
\end{theorem}

\vspace*{-2mm}
\subsection{Annihilation Operator Only Systems}
Annihilation operator only linear quantum systems are a special case of the above class of linear quantum systems, where the QSDEs (\ref{eq:lqs_1}) can be described purely in terms of the vector of annihilation operators \cite{MP1,MP2}:
\vspace*{-6mm}\small
\begin{equation}\label{eq:ann_lqs_1}
\begin{split}
da(t) &= F_1a(t)dt + G_1d\mathcal{A}(t);\\
d\mathcal{A}^{out}(t) &= H_1a(t)dt + K_1d\mathcal{A}(t).
\end{split}
\end{equation}
\normalsize \vspace*{-4mm}

\vspace*{-2mm}
\begin{theorem}\label{thm:ann_phys_rlz}
(See \cite{MP1,IRP3}) An annihilation operator only linear quantum system of the form (\ref{eq:ann_lqs_1}) is physically realizable, if and only if there exists a complex commutation matrix $\Theta = \Theta^\dagger > 0$, satisfying
\vspace*{-6mm}\small
\begin{equation}\label{eq:comm_rel2}
\Theta = \left[a,a^\dagger \right],
\end{equation}
\normalsize \vspace*{-6mm}

such that
\vspace*{-1.2cm}\small
\begin{equation}
\begin{split}
F_1\Theta + \Theta F_1^\dagger + G_1G_1^\dagger &= 0,\\
G_1 &= -\Theta H_1^\dagger ,\\
K_1 &= I.
\end{split}
\end{equation}
\normalsize \vspace*{-5mm}
\end{theorem}

\vspace*{-2mm}
\subsection{Linear Quantum System from Quantum Optics}
\begin{figure}
\centering
\includegraphics[width=0.3\textwidth]{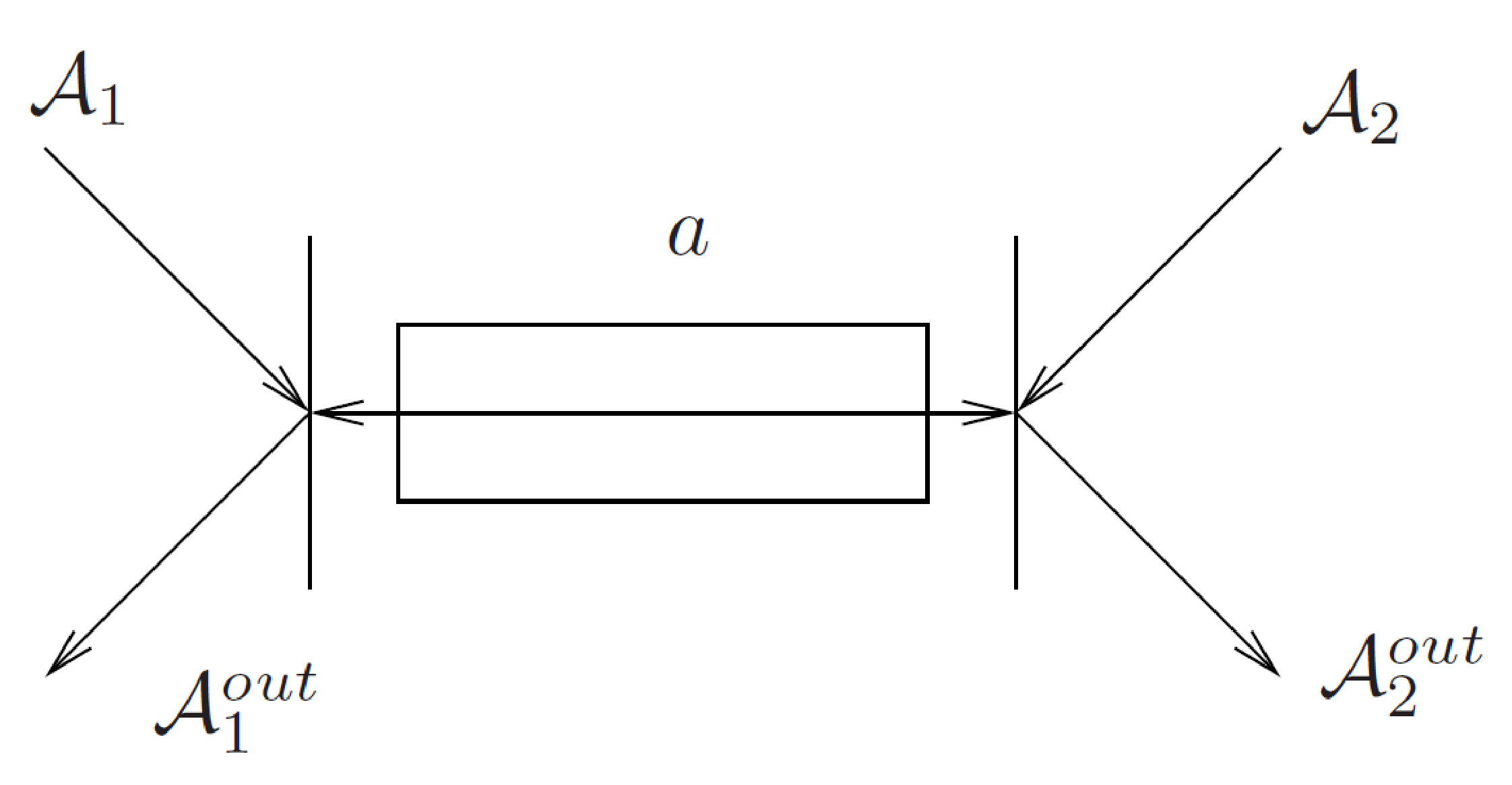}
\caption{\small Schematic diagram of a dynamic optical squeezer.}
\label{fig:sqz_scm}
\end{figure}

An example of a linear quantum system is a linearized dynamic optical squeezer. This is an optical cavity with a non-linear optical element inside as shown in Fig. \ref{fig:sqz_scm}. Such a dynamic squeezer can be described by the quantum stochastic differential equations \cite{IRP2}:
\vspace*{-6mm}\small
\begin{equation}
\begin{split}
da &= -\frac{\gamma}{2}adt - \chi a^{*} dt - \sqrt{\kappa_1}d\mathcal{A}_1 - \sqrt{\kappa_2}d\mathcal{A}_2;\\
d\mathcal{A}_1^{out} &= \sqrt{\kappa_1}adt + d\mathcal{A}_1;\\
d\mathcal{A}_2^{out} &= \sqrt{\kappa_2}adt + d\mathcal{A}_2,
\end{split}
\end{equation}
\normalsize \vspace*{-4mm}

where $\kappa_1,\kappa_2>0$, $\chi \in \mathbb{C}$, and $a$ is a single annihilation operator of the cavity mode \cite{BR,GZ}. This leads to a linear quantum system of the form (\ref{eq:lqs_1}) as follows:
\vspace*{-8mm}\small
\begin{equation}\label{eq:ldos_qo}
\begin{split}
\left[\begin{array}{c}
da(t)\\
da(t)^{*}
\end{array}\right] &= \left[\begin{array}{cc}
-\frac{\gamma}{2} & -\chi\\
-\chi^{*} & -\frac{\gamma}{2}
\end{array}\right] \left[\begin{array}{c}
a(t)\\
a(t)^{*}
\end{array}\right] dt\\
&- \sqrt{\kappa_1} \left[\begin{array}{c}
d\mathcal{A}_1(t)\\
d\mathcal{A}_1(t)^{*}
\end{array}\right] - \sqrt{\kappa_2} \left[\begin{array}{c}
d\mathcal{A}_2(t)\\
d\mathcal{A}_2(t)^{*}
\end{array}\right];\\
\left[\begin{array}{c}
d\mathcal{A}_1^{out}(t)\\
d\mathcal{A}_1^{out}(t)^{*}
\end{array}\right] &= \sqrt{\kappa_1} \left[\begin{array}{c}
a(t)\\
a(t)^{*}
\end{array}\right] dt + \left[\begin{array}{c}
d\mathcal{A}_1(t)\\
d\mathcal{A}_1(t)^{*}
\end{array}\right];\\
\left[\begin{array}{c}
d\mathcal{A}_2^{out}(t)\\
d\mathcal{A}_2^{out}(t)^{*}
\end{array}\right] &= \sqrt{\kappa_2} \left[\begin{array}{c}
a(t)\\
a(t)^{*}
\end{array}\right] dt + \left[\begin{array}{c}
d\mathcal{A}_2(t)\\
d\mathcal{A}_2(t)^{*}
\end{array}\right].
\end{split}
\end{equation}
\normalsize \vspace*{-6mm}

The above quantum system requires $\gamma = \kappa_1 + \kappa_2$ in order for the system to be physically realizable.

Also, the above quantum optical system can be described purely in terms of the annihilation operator, if and only if $\chi = 0$, i.e. there is no squeezing, in which case it reduces to a passive optical cavity. This leads to a linear quantum system of the form (\ref{eq:ann_lqs_1}) as follows:
\vspace*{-6mm}\small
\begin{equation}
\begin{split}
da &= -\frac{\gamma}{2}adt - \sqrt{\kappa_1} d\mathcal{A}_1 - \sqrt{\kappa_2} d\mathcal{A}_2;\\
d\mathcal{A}_1^{out} &= \sqrt{\kappa_1} adt + d\mathcal{A}_1;\\
d\mathcal{A}_2^{out} &= \sqrt{\kappa_2} adt + d\mathcal{A}_2,
\end{split}
\end{equation}
\normalsize \vspace*{-6mm}

where again the system is physically realizable when we have $\gamma = \kappa_1 + \kappa_2$.

\section{Purely-Classical Estimation}\label{sec:pce}
\begin{figure}[!b]
\centering
\includegraphics[width=0.5\textwidth]{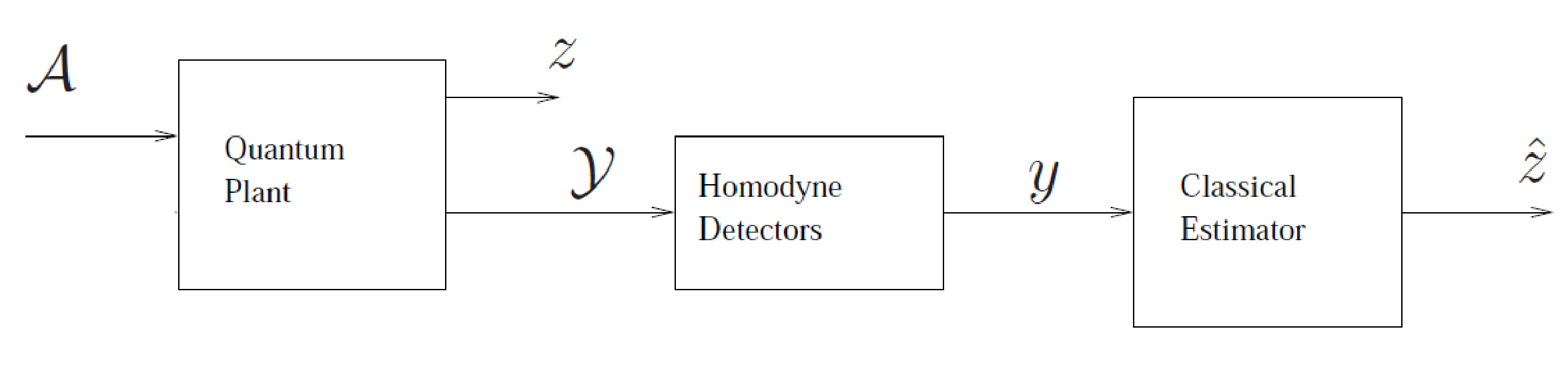}
\caption{\small Schematic diagram of purely-classical estimation.}
\label{fig:cls_scm}
\end{figure}

\vspace*{-3mm}
The schematic diagram of a purely-classical estimation scheme is provided in Fig. \ref{fig:cls_scm}. We consider a quantum plant of the form (\ref{eq:lqs_1}), (\ref{eq:lqs_2}), defined as follows:
\vspace*{-8mm}\small
\begin{equation}\label{eq:plant}
\begin{split}
\left[\begin{array}{c}
da\\
da^{\#}
\end{array}\right] &= F \left[\begin{array}{c}
a\\
a^{\#}
\end{array}\right] dt + G \left[\begin{array}{c}
d\mathcal{A}\\
d\mathcal{A}^{\#} 
\end{array}\right];\\
\left[\begin{array}{c}
d\mathcal{Y}\\
d\mathcal{Y}^{\#}
\end{array}\right] &= H \left[\begin{array}{c}
a\\
a^{\#}
\end{array}\right] dt + K \left[\begin{array}{c}
d\mathcal{A}\\
d\mathcal{A}^{\#}
\end{array}\right];\\
z &= C\left[\begin{array}{c}
a\\
a^{\#} 
\end{array}\right].
\end{split}
\end{equation}
\normalsize \vspace*{-6mm}

Here, $z$ denotes a scalar operator on the underlying Hilbert space and represents the quantity to be estimated. Also, $\mathcal{Y}$ is the vector of output fields of the plant, and $\mathcal{A}$ is a vector of quantum noises acting on the plant.

In purely-classical estimation, a quadrature of each component of the vector $\mathcal{Y}$ is measured using homodyne detection to produce a corresponding classical signal $y_i$:
\vspace*{-7mm}\small
\begin{equation}\label{eq:class_hd}
\begin{split}
dy_1 &= \frac{1}{2}e^{-\iota\theta_1}d\mathcal{Y}_1 + \frac{1}{2}e^{\iota\theta_1}d\mathcal{Y}_1^{*};\\
&\vdots\\
dy_m &= \frac{1}{2}e^{-\iota\theta_m}d\mathcal{Y}_m + \frac{1}{2}e^{\iota\theta_m}d\mathcal{Y}_m^{*}.
\end{split}
\end{equation}
\normalsize \vspace*{-6mm}

Here, $\iota=\sqrt{-1}$, and the angles $\theta_1,\hdots,\theta_m$ determine the quadrature measured by each homodyne detector. The vector of real classical signals $y = [y_1 \hdots y_m]^T$ is then used as the input to a classical estimator defined as:
\vspace*{-6mm}\small
\begin{equation}\label{eq:class_estimator}
\begin{split}
dx_e &= F_ex_edt + G_edy;\\
\hat{z} &= H_ex_e.
\end{split}
\end{equation}
\normalsize \vspace*{-6mm}

Here $\hat{z}$ is a scalar classical estimate of the quantity $z$. The estimation error corresponding to this estimate is

\vspace*{-6mm}\small
\begin{equation}\label{eq:est_err}
e = z - \hat{z}.
\end{equation}
\normalsize \vspace*{-6mm}

Then, the optimal classical estimator is defined as the system (\ref{eq:class_estimator}) that minimizes the quantity
\vspace*{-6mm}\small
\begin{equation}\label{eq:est_cost}
\bar{J}_c = \lim_{t \to \infty} \left\langle e^{*}(t)e(t) \right\rangle ,
\end{equation}
\normalsize \vspace*{-6mm}

which is the mean-square error of the estimate. Here, $\langle \cdot \rangle$ denotes the quantum expectation over the joint quantum-classical system defined by (\ref{eq:plant}), (\ref{eq:class_hd}), (\ref{eq:class_estimator}).

The optimal classical estimator is given by the standard (complex) Kalman filter defined for the system (\ref{eq:plant}), (\ref{eq:class_hd}). This optimal classical estimator is obtained from the solution to the algebraic Riccati equation:

\vspace*{-8mm}\small
\begin{equation}\label{eq:riccati}
\begin{split}
F_a\bar{P}_e &+ \bar{P}_eF_a^\dagger + G_aG_a^\dagger - (G_aK_a^\dagger + \bar{P}_eH_a^\dagger)L^\dagger \\
&\times(LK_aK_a^\dagger L^\dagger)^{-1}L(G_aK_a^\dagger + \bar{P}_eH_a^\dagger)^\dagger = 0,
\end{split}
\end{equation}
\normalsize \vspace*{-6mm}
where
\vspace*{-7mm}\small
\begin{equation}
\begin{split}
\quad \, \, \, \, \, \quad F_a &= F, \quad G_a = G,\\
H_a &= H, \quad K_a = K, \quad L = \left[\begin{array}{cc}
L_1 & L_2
\end{array}\right],\\
L_1 &= \left[\begin{array}{cccc}
\frac{1}{2}e^{-\iota\theta_1} & 0 & \hdots & 0\\
0 & \frac{1}{2}e^{-\iota\theta_2} & \hdots & 0\\
0 & 0 & \ddots & 0\\
0 & \hdots & 0 & \frac{1}{2}e^{-\iota\theta_m}
\end{array}\right],\\
L_2 &= \left[\begin{array}{cccc}
\frac{1}{2}e^{\iota\theta_1} & 0 & \hdots & 0\\
0 & \frac{1}{2}e^{\iota\theta_2} & \hdots & 0\\
0 & 0 & \ddots & 0\\
0 & \hdots & 0 & \frac{1}{2}e^{\iota\theta_m}
\end{array}\right].
\end{split}
\end{equation}
\normalsize \vspace*{-6mm}

Here we have assumed that the quantum noise $\mathcal{A}$ is purely canonical, i.e. $d\mathcal{A}d\mathcal{A}^\dagger = Idt$ and hence $K=I$.

Equation (\ref{eq:riccati}) thus becomes:
\vspace*{-6mm}\small
\begin{equation}\label{eq:class_riccati}
\begin{split}
F\bar{P}_e &+ \bar{P}_eF^\dagger + GG^\dagger - (G + \bar{P}_eH^\dagger)L^\dagger \\
&\times (LL^\dagger)^{-1}L(G + \bar{P}_eH^\dagger)^\dagger =0.
\end{split}
\end{equation}
\normalsize \vspace*{-7mm}

Then, the corresponding optimal classical estimator (\ref{eq:class_estimator}) is defined by the equations:

\vspace*{-8mm}\small
\begin{equation}\label{eq:cls_sys_est}
\begin{split}
F_e &= F - G_eLH;\\
G_e &= (G + \bar{P}_eH^\dagger)L^\dagger (LL^\dagger)^{-1};\\
H_e &= C.
\end{split}
\end{equation}
\normalsize \vspace*{-7mm}

The value of the cost (\ref{eq:est_cost}) is given by
\vspace*{-6mm}\small
\begin{equation}\label{eq:class_cost}
\bar{J}_c = C\bar{P}_eC^\dagger .
\end{equation}
\normalsize \vspace*{-6mm}

\section{Coherent-Classical Estimation}\label{sec:cce}
\begin{figure}
\centering
\includegraphics[width=0.5\textwidth]{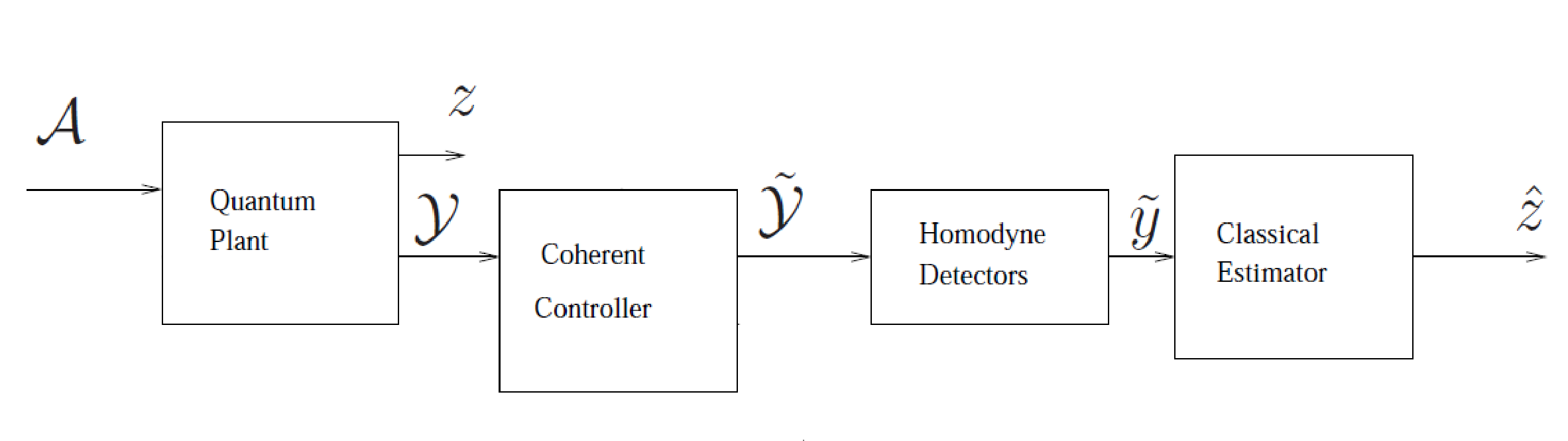}
\caption{\small Schematic diagram of coherent-classical estimation.}
\label{fig:coh_cls_scm}
\end{figure}

\vspace*{-3mm}
In coherent-classical estimation scheme of Fig. \ref{fig:coh_cls_scm}, the plant output $\mathcal{Y}(t)$ does not directly drive a bank of homodyne detectors as in (\ref{eq:class_hd}). Rather, this is fed into another quantum system called a coherent controller, defined as
\vspace*{-8mm}\small
\begin{equation}\label{eq:coh_controller}
\begin{split}
\left[\begin{array}{c}
da_c\\
da_c^{\#}
\end{array}\right] &= F_c \left[\begin{array}{c}
a_c\\
a_c^{\#}
\end{array}\right] dt + G_c \left[\begin{array}{c}
d\mathcal{Y}\\
d\mathcal{Y}^{\#} 
\end{array}\right];\\
\left[\begin{array}{c}
d\tilde{\mathcal{Y}}\\
d\tilde{\mathcal{Y}}^{\#}
\end{array}\right] &= H_c \left[\begin{array}{c}
a_c\\
a_c^{\#}
\end{array}\right] dt + K_c \left[\begin{array}{c}
d\mathcal{Y}\\
d\mathcal{Y}^{\#} 
\end{array}\right].
\end{split}
\end{equation}
\normalsize \vspace*{-6mm}

A quadrature of each component of $\tilde{\mathcal{Y}}$ is homodyne detected to yield a corresponding classical signal $\tilde{y}_i$:
\vspace*{-6mm}\small
\begin{equation}\label{eq:coh_class_hd}
\begin{split}
d\tilde{y}_1 &= \frac{1}{2}e^{-\iota\tilde{\theta}_1}d\tilde{\mathcal{Y}}_1 + \frac{1}{2}e^{\iota\tilde{\theta}_1}d\tilde{\mathcal{Y}}_1^{*};\\
&\vdots\\
d\tilde{y}_{\tilde{m}} &= \frac{1}{2}e^{-\iota\tilde{\theta}_{\tilde{m}}}d\tilde{\mathcal{Y}}_{\tilde{m}} + \frac{1}{2}e^{\iota\tilde{\theta}_{\tilde{m}}}d\tilde{\mathcal{Y}}_{\tilde{m}}^{*}.
\end{split}
\end{equation}
\normalsize \vspace*{-4mm}

Here, the angles $\tilde{\theta}_1,\hdots,\tilde{\theta}_{\tilde{m}}$ determine the quadrature measured by each homodyne detector. The vector of real classical signals $\tilde{y} = [\tilde{y}_1 \hdots \tilde{y}_{\tilde{m}}]^T$ is then used as the input to a classical estimator defined as follows:
\vspace*{-6mm}\small
\begin{equation}\label{eq:coh_class_estimator}
\begin{split}
d\tilde{x}_e &= \tilde{F}_e\tilde{x}_edt + \tilde{G}_ed\tilde{y};\\
\hat{z} &= \tilde{H}_e\tilde{x}_e.
\end{split}
\end{equation}
\normalsize \vspace*{-4mm}

Here $\hat{z}$ is a scalar classical estimate of the quantity $z$. Corresponding to this estimate is the estimation error (\ref{eq:est_err}). Then, the optimal coherent-classical estimator is defined as the systems (\ref{eq:coh_controller}), (\ref{eq:coh_class_estimator}) which together minimize the quantity (\ref{eq:est_cost}). Note that the coherent controller does not directly produce an estimate of a plant variable as in the quantum observer of Ref. \cite{MJ}. Instead, it only works in combination with the classical estimator to yield a classical estimate of the quantity $z$.

We can now combine the quantum plant (\ref{eq:plant}) and the coherent controller (\ref{eq:coh_controller}) to yield an augmented quantum linear system defined by the following QSDEs:
\vspace*{-8mm}\small
\begin{equation}\label{eq:coh_system}
\begin{split}
\left[\begin{array}{c}
da\\
da^{\#}\\
da_c\\
da_c^{\#}
\end{array}\right] &= \left[\begin{array}{cc}
F & 0\\
G_cH & F_c
\end{array}\right] \left[\begin{array}{c}
a\\
a^{\#}\\
a_c\\
a_c^{\#}
\end{array}\right] dt+ \left[\begin{array}{c}
G\\
G_cK
\end{array}\right] \left[\begin{array}{c}
d\mathcal{A}\\
d\mathcal{A}^{\#}
\end{array}\right];\\
\left[\begin{array}{c}
d\tilde{\mathcal{Y}}\\
d\tilde{\mathcal{Y}}^{\#}
\end{array}\right] &= \left[\begin{array}{cc}
K_cH & H_c
\end{array}\right] \left[\begin{array}{c}
a\\
a^{\#}\\
a_c\\
a_c^{\#}
\end{array}\right] dt + K_cK \left[\begin{array}{c}
d\mathcal{A}\\
d\mathcal{A}^{\#}
\end{array}\right].
\end{split}
\end{equation}
\normalsize \vspace*{-6mm}

The optimal classical estimator is given by the standard (complex) Kalman filter defined for the system (\ref{eq:coh_system}), (\ref{eq:coh_class_hd}). This optimal classical estimator is obtained from the solution $\tilde{P}_e$ to an algebraic Riccati equation of the form (\ref{eq:riccati}), where
\vspace*{-8mm}\small
\begin{equation}\label{eq:coh_sys_matrices}
\begin{split}
F_a &= \left[\begin{array}{cc}
F & 0\\
G_cH & F_c
\end{array}\right], \quad G_a = \left[\begin{array}{c}
G\\
G_cK
\end{array}\right],\\
H_a &= \left[\begin{array}{cc}
K_cH & H_c
\end{array}\right], \quad K_a = K_cK, \quad L = \left[\begin{array}{cc}
\tilde{L}_1 & \tilde{L}_2
\end{array}\right],\\
\tilde{L}_1 &= \left[\begin{array}{cccc}
\frac{1}{2}e^{-\iota\tilde{\theta}_1} & 0 & \hdots & 0\\
0 & \frac{1}{2}e^{-\iota\tilde{\theta}_2} & \hdots & 0\\
0 & 0 & \ddots & 0\\
0 & \hdots & 0 & \frac{1}{2}e^{-\iota\tilde{\theta}_{\tilde{m}}}
\end{array}\right],\\
\tilde{L}_2 &= \left[\begin{array}{cccc}
\frac{1}{2}e^{\iota\tilde{\theta}_1} & 0 & \hdots & 0\\
0 & \frac{1}{2}e^{\iota\tilde{\theta}_2} & \hdots & 0\\
0 & 0 & \ddots & 0\\
0 & \hdots & 0 & \frac{1}{2}e^{\iota\tilde{\theta}_{\tilde{m}}}
\end{array}\right],
\end{split}
\end{equation}
\normalsize \vspace*{-7mm}

where since the quantum noise $\mathcal{A}$ is assumed to be purely canonical, i.e. $d\mathcal{A}d\mathcal{A}^\dagger = Idt$, we have $K_a=K_cK=I$, which requires $K_c=I$ too, as $K=I$.

Then, the corresponding optimal classical estimator (\ref{eq:coh_class_estimator}) is defined by the equations:
\vspace*{-8mm}\small
\begin{equation}\label{eq:coh_sys_est}
\begin{split}
\tilde{F}_e &= F_a - \tilde{G}_eLH_a;\\
\tilde{G}_e &= (G_aK_a^\dagger + \tilde{P}_eH_a^\dagger)L^\dagger(LL^\dagger)^{-1};\\
\tilde{H}_e &= \left[\begin{array}{cc}
C & 0
\end{array}\right].
\end{split}
\end{equation}
\normalsize \vspace*{-6mm}

We write:
\vspace*{-1cm}\small
\begin{equation}\label{eq:coh_err_cov}
\tilde{P}_e = \left[\begin{array}{cc}
P_1 & P_2\\
P_2^\dagger & P_3
\end{array}\right],
\end{equation}
\normalsize \vspace*{-6mm}

where $P_1$ is of the same dimension as $\bar{P}_e$. Then, the corresponding cost of the form (\ref{eq:est_cost}) is
\vspace*{-6mm}\small
\begin{equation}\label{eq:coh_cost}
\tilde{J}_c = \left[\begin{array}{cc}
C & 0
\end{array}\right] \tilde{P}_e \left[\begin{array}{c}
C^\dagger \\
0
\end{array}\right] = CP_1C^\dagger .
\end{equation}
\normalsize \vspace*{-6mm}

Thus, the optimal coherent-classical estimation problem can be solved by first choosing the coherent controller (\ref{eq:coh_controller}) to minimize the cost (\ref{eq:coh_cost}). Then, the classical estimator (\ref{eq:coh_class_estimator}) is constructed according to (\ref{eq:coh_sys_est}).

\begin{remark}
Note that the combined plant-controller system being measured here is a fully quantum system, such that the controller preserves the quantum coherence of the quantum plant output (that is not measured directly) and yet can be chosen suitably, as mentioned above, to assist in improving the precision of the classical estimate of a plant variable. On the other hand, with purely-classical estimation, we have no control over the variables of the quantum system (just the plant itself) being measured.
\end{remark}

\begin{theorem}\label{thm:central_result1}
Consider a coherent-classical estimation scheme defined by (\ref{eq:ann_lqs_1}) ($\mathcal{A}^{out}$ being $\mathcal{Y}$), (\ref{eq:coh_controller}), (\ref{eq:coh_class_hd}) and (\ref{eq:coh_class_estimator}), such that the plant is physically realizable, with the cost $\tilde{J}_c$ defined in (\ref{eq:coh_cost}). Also, consider the corresponding purely-classical estimation scheme defined by (\ref{eq:ann_lqs_1}), (\ref{eq:class_hd}) and (\ref{eq:class_estimator}), such that the plant is physically realizable, with the cost $\bar{J}_c$ defined in (\ref{eq:class_cost}). Then,
\vspace*{-6mm}\small
\begin{equation}
\tilde{J}_c = \bar{J}_c.
\end{equation}
\normalsize \vspace*{-4mm}
\end{theorem}
\vspace*{-2mm}
\begin{proof}
We first consider the form of the system (\ref{eq:plant}) with the assumption that the plant is an annihilation operator only system. A quantum system (\ref{eq:lqs_1}), (\ref{eq:lqs_2}) is characterized by annihilation operators only when $F_2, G_2, H_2, K_2 = 0$.

Then, the equations for the annihilation operators in (\ref{eq:plant}) take the form
\vspace*{-7mm}\small
\begin{equation}\label{eq:ann_plant}
\begin{split}
da &= F_1adt + G_1d\mathcal{A};\\
d\mathcal{Y} &= H_1adt + K_1d\mathcal{A}.
\end{split}
\end{equation}
\normalsize \vspace*{-5mm}

The corresponding equations for creation operators are
\vspace*{-7mm}\small
\begin{equation}\label{eq:crn_plant}
\begin{split}
da^{\#} &= F_1^{\#}a^{\#}dt + G_1^{\#}d\mathcal{A}^{\#};\\
d\mathcal{Y}^{\#} &= H_1^{\#}a^{\#}dt + K_1^{\#}d\mathcal{A}^{\#},
\end{split}
\end{equation}
\normalsize \vspace*{-4mm}
Hence, the plant is described by (\ref{eq:plant}), where
\vspace*{-6mm}\small
\begin{equation}
\begin{split}
F &= \left[\begin{array}{cc}
F_1 & 0\\
0 & F_1^{\#}
\end{array}\right], G = \left[\begin{array}{cc}
G_1 & 0\\
0 & G_1^{\#}
\end{array}\right],\\
H &= \left[\begin{array}{cc}
H_1 & 0\\
0 & H_1^{\#}
\end{array}\right], K = \left[\begin{array}{cc}
K_1 & 0\\
0 & K_1^{\#}
\end{array}\right].
\end{split}
\end{equation}
\normalsize \vspace*{-5mm}

Next, we use the assumption that the plant is physically realizable. Then, by applying Theorem \ref{thm:ann_phys_rlz} to (\ref{eq:ann_plant}), there exists a matrix $\Theta_1>0$, such that
\vspace*{-7mm}\small
\begin{equation}\label{eq:ann_phys_rlz}
\begin{split}
F_1\Theta_1 +\Theta_1 F_1^\dagger +G_1G_1^\dagger &=0,\\
G_1&=-\Theta_1 H_1^\dagger ,\\
K_1&=I.
\end{split}
\end{equation}
\normalsize \vspace*{-4mm}
Hence,
\vspace*{-1cm}\small
\begin{equation}\label{eq:crn_phys_rlz}
\begin{split}
F_1^{\#}\Theta_1^{\#} +\Theta_1^{\#} F_1^T +G_1^{\#}G_1^T &=0,\\
G_1^{\#}&=-\Theta_1^{\#} H_1^T ,\\
K_1^{\#}&=I.
\end{split}
\end{equation}
\normalsize \vspace*{-4mm}

Combining (\ref{eq:ann_phys_rlz}) and (\ref{eq:crn_phys_rlz}), we get
\vspace*{-7mm}\small
\begin{equation}\label{eq:plant_phys_rlz}
\begin{split}
F\Theta +\Theta F^\dagger +GG^\dagger &=0,\\
G&=-\Theta H^\dagger ,\\
K&=I,
\end{split}
\end{equation}
\normalsize \vspace*{-4mm}
where $\Theta = \left[\begin{array}{cc}
\Theta_1 & 0\\
0 & \Theta_1^{\#}
\end{array}\right] > 0$. Clearly, $\bar{P}_e = \Theta$ satisfies (\ref{eq:class_riccati}) owing to (\ref{eq:plant_phys_rlz}). Also, $\tilde{P}_e = \left[\begin{array}{cc}
\Theta & 0\\
0 & P_3 \end{array}\right]$ satisfies (\ref{eq:riccati}), (\ref{eq:coh_sys_matrices}) for the coherent-classical estimation case. Here, $P_3>0$ is the error-covariance of the purely-classical estimation of the coherent controller alone. Thus, we get $\bar{J}_c = \tilde{J}_c = C\Theta C^\dagger $.
\hfill \qed
\end{proof}

\begin{remark}
Note that the Kalman gain of the purely-classical estimator is $0$ when $\bar{P}_e = \Theta$, i.e. the Kalman state estimate is independent of the measurement. This is consistent with Cor. 1 of Ref. \cite{IRP3}, which states that for a physically realizable annihilation operator quantum system with only quantum noise inputs, any output contains no information about the system's internal variables.
\end{remark}

\begin{remark}
Theorem \ref{thm:central_result1} implies that coherent-classical estimation of a physically realizable annihilation operator quantum plant performs identical to, and no better than, purely-classical estimation of the plant. This is so because the output field of the plant contains no information about the plant's internal variables and, thus, simply serves as a quantum white noise input for the controller.
\end{remark}

Now, we present an example to illustrate Theorem \ref{thm:central_result1}. Let the quantum plant be a dynamic squeezer (See (\ref{eq:ldos_qo})):
\vspace*{-7mm}\small
\begin{equation}\label{eq:sqz_plant}
\begin{split}
\left[\begin{array}{c}
da\\
da^{*}
\end{array}\right] &= \left[\begin{array}{cc}
-\frac{\gamma}{2} & -\chi\\
-\chi^{*} & -\frac{\gamma}{2}
\end{array}\right] \left[\begin{array}{c}
a\\
a^{*}
\end{array}\right] dt - \sqrt{\kappa} \left[\begin{array}{c}
d\mathcal{A}\\
d\mathcal{A}^{*}
\end{array}\right];\\
\left[\begin{array}{c}
d\mathcal{Y}\\
d\mathcal{Y}^{*}
\end{array}\right] &= \sqrt{\kappa} \left[\begin{array}{c}
a\\
a^{*}
\end{array}\right] dt + \left[\begin{array}{c}
d\mathcal{A}\\
d\mathcal{A}^{*}
\end{array}\right];\\
z &= \left[\begin{array}{cc}
0.2 & -0.2
\end{array}\right] \left[\begin{array}{c}
a\\
a^{*}
\end{array}\right].
\end{split}
\end{equation}
\normalsize \vspace*{-5mm}

Here, we choose $\gamma = 4$, $\kappa = 4$ and $\chi = 0$. Note that this system is physically realizable, since $\gamma = \kappa$, and is annihilation operator only, since $\chi = 0$. In fact, this system corresponds to a passive optical cavity. We then calculate the optimal classical-only state estimator and the error $\bar{J}_c$ of (\ref{eq:class_cost}) for this system using the standard Kalman filter equations corresponding to homodyne detector angles varying from $\theta = 0^{\circ}$ to $\theta = 180^{\circ}$.\vspace*{-2mm}

We next consider coherent-classical estimation, where the coherent controller (\ref{eq:coh_controller}) is also a dynamic squeezer:

\vspace*{-8mm}\small
\begin{equation}\label{eq:sqz_ctrlr}
\begin{split}
\left[\begin{array}{c}
da\\
da^{*}
\end{array}\right] &= \left[\begin{array}{cc}
-\frac{\gamma}{2} & -\chi\\
-\chi^{*} & -\frac{\gamma}{2}
\end{array}\right] \left[\begin{array}{c}
a\\
a^{*}
\end{array}\right] dt - \sqrt{\kappa} \left[\begin{array}{c}
d\mathcal{Y}\\
d\mathcal{Y}^{*}
\end{array}\right];\\
\left[\begin{array}{c}
d\tilde{\mathcal{Y}}\\
d\tilde{\mathcal{Y}}^{*}
\end{array}\right] &= \sqrt{\kappa} \left[\begin{array}{c}
a\\
a^{*}
\end{array}\right] dt + \left[\begin{array}{c}
d\mathcal{Y}\\
d\mathcal{Y}^{*}
\end{array}\right].
\end{split}
\end{equation}
\normalsize \vspace*{-5mm}

\begin{figure}
\centering
\includegraphics[width=0.5\textwidth]{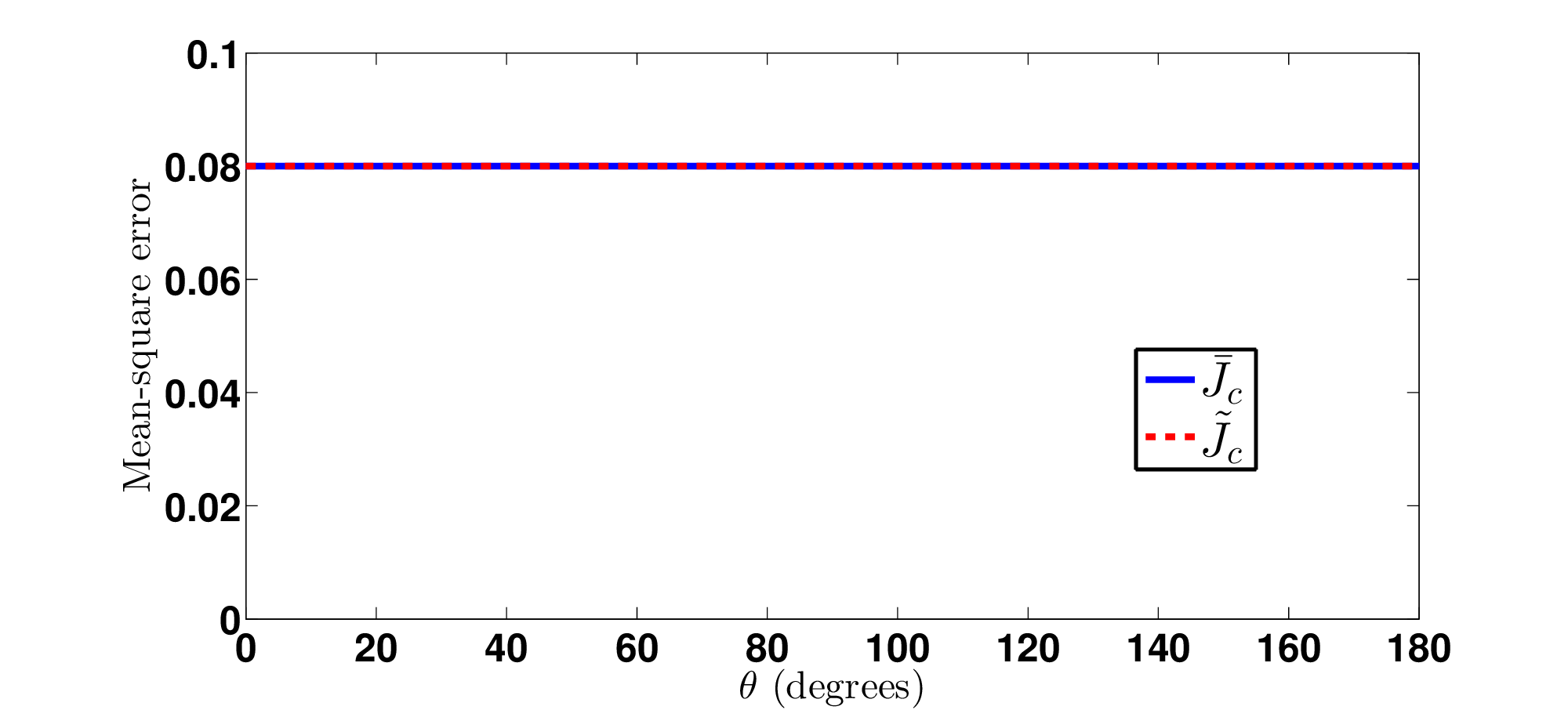}
\caption{\small Estimation error vs. homodyne angle $\theta$ in the case of an annihilation operator only plant.}
\label{fig:cls_vs_coh_cls_1}
\end{figure}

Here, we choose $\gamma = 16$, $\kappa = 16$ and $\chi = 2$, so that the system is physically realizable. Then, the classical estimator for this case is calculated according to (\ref{eq:coh_system}), (\ref{eq:coh_sys_matrices}), (\ref{eq:riccati}), (\ref{eq:coh_sys_est}) for the homodyne detector angle varying from $\theta=0^{\circ}$ to $\theta=180^{\circ}$. The resulting cost $\tilde{J}_c$ in (\ref{eq:coh_cost}) alongwith the cost for the purely-classical estimator is shown in Fig. \ref{fig:cls_vs_coh_cls_1}. Clearly, both the classical-only and coherent-classical estimators have the same estimation error cost for all homodyne angles. This illustrates Theorem \ref{thm:central_result1}.\vspace*{-2mm}

\begin{figure}[!b]
\centering
\includegraphics[width=0.5\textwidth]{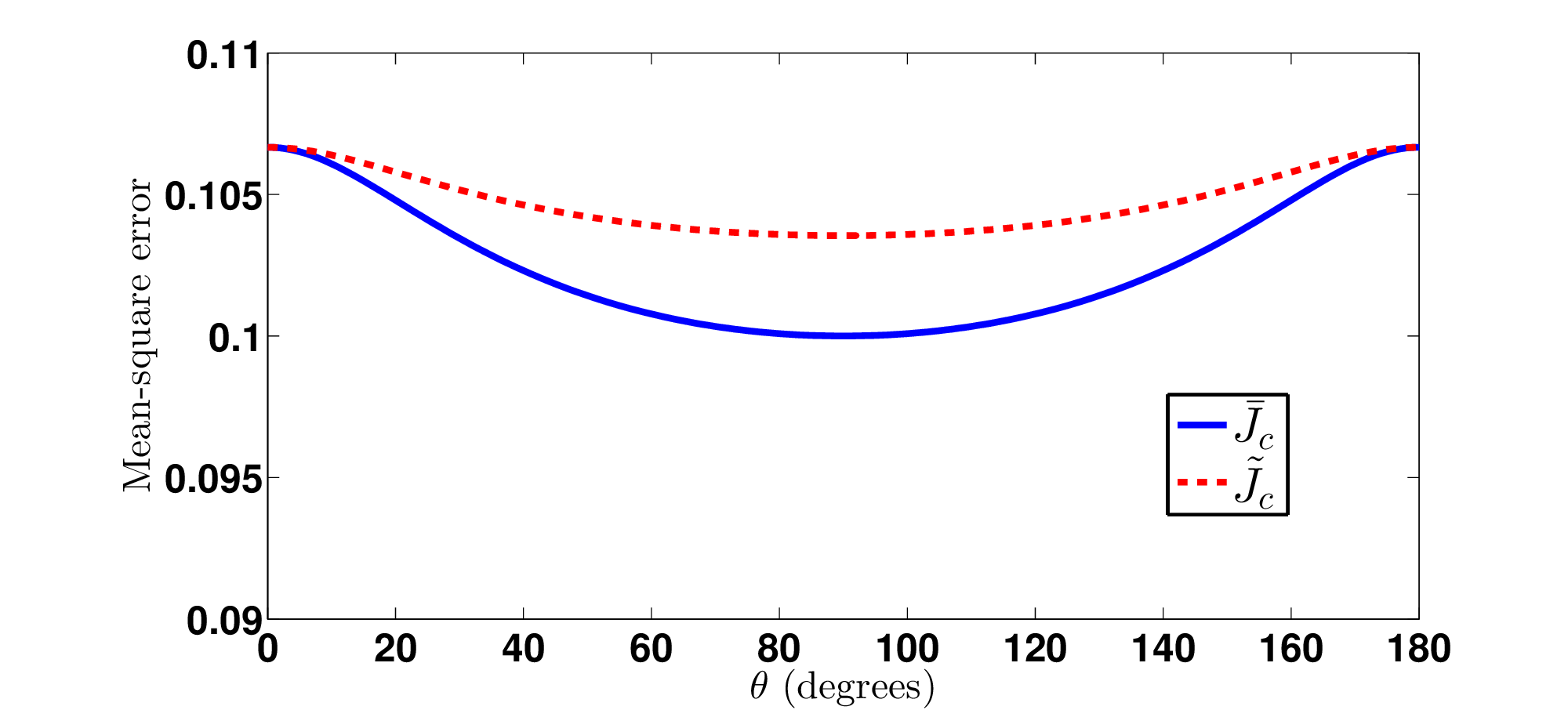}
\caption{\small Estimation error vs. homodyne angle $\theta$ in the case of an annihilation operator only controller.}
\label{fig:cls_vs_coh_cls_2}
\end{figure}

Next, we consider a case where the controller is a physically realizable annihilation operator only system. But, the plant is physically realizable and has $\chi \neq 0$. In (\ref{eq:sqz_plant}), we choose $\gamma = 4$, $\kappa = 4$, $\chi = 0.5$, and in (\ref{eq:sqz_ctrlr}), $\gamma = 16$, $\kappa = 16$, $\chi = 0$. Fig. \ref{fig:cls_vs_coh_cls_2} then shows that the coherent-classical error is greater than or equal to the purely-classical error for all homodyne angles.\vspace*{-2mm}

In fact, we observe that the coherent-classical estimator can be no better than the purely-classical estimator, when the coherent controller is an annihilation operator only system. We present this here as a conjecture, which is a consequence of the quantum data processing inequality from Ref. \cite{SN}. Indeed, the coherent information in the plant cannot be increased by additional dynamics of a coherent controller, that does not have any squeezing and provides no feedback to the plant. Thus, such a controller cannot improve the estimation accuracy.\vspace*{-2mm}

\begin{conjecture}\label{thm:central_result2}
Consider a coherent-classical estimation scheme defined by (\ref{eq:plant}), (\ref{eq:coh_controller}), (\ref{eq:coh_class_hd}) and (\ref{eq:coh_class_estimator}), where the plant is physically realizable and the coherent controller is a physically realizable annihilation operator only system, with the cost $\tilde{J}_c$ defined in (\ref{eq:coh_cost}). Also, consider the corresponding purely-classical estimation scheme defined by (\ref{eq:plant}), (\ref{eq:class_hd}) and (\ref{eq:class_estimator}), such that the plant is physically realizable, with the cost $\bar{J}_c$ defined in (\ref{eq:class_cost}). Then,
\vspace*{-6mm}\small
\begin{equation}
\tilde{J}_c \geq \bar{J}_c.
\end{equation}
\normalsize \vspace*{-6mm}
\end{conjecture}

\begin{figure}
\centering
\includegraphics[width=0.5\textwidth]{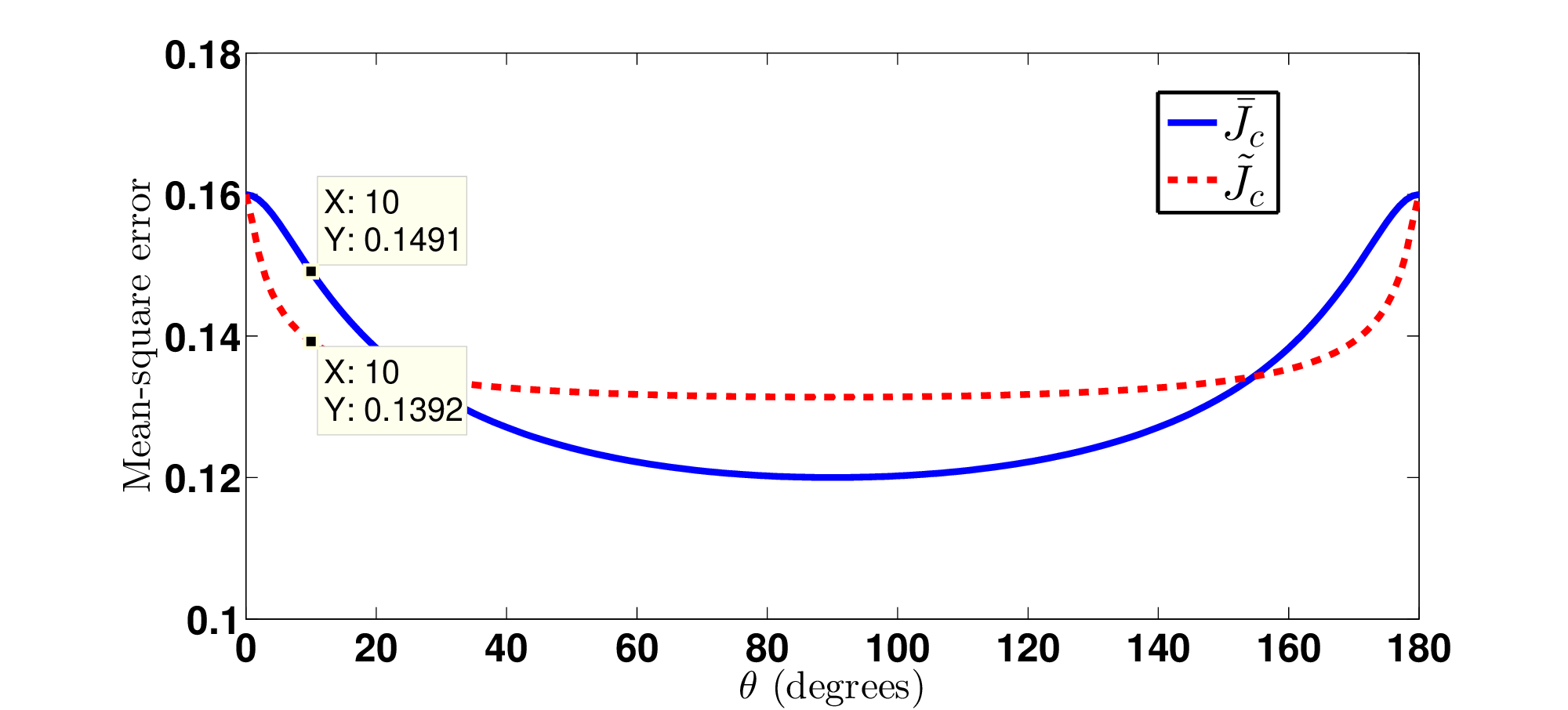}
\caption{\small Estimation error vs. homodyne angle $\theta$ in the case of a squeezer plant and a squeezer controller.}
\label{fig:cls_vs_coh_cls_3}
\end{figure}

Furthermore, we see an example where both the plant and the controller are physically realizable quantum systems with $\chi \neq 0$. In (\ref{eq:sqz_plant}), we choose $\gamma = 4$, $\kappa = 4$, $\chi = 1$, and in (\ref{eq:sqz_ctrlr}), $\gamma = 16$, $\kappa = 16$, $\chi = 4$. Fig.~\ref{fig:cls_vs_coh_cls_3} then shows that the coherent-classical estimator can perform better than the purely-classical estimator, e.g.~for a homodyne angle of $\theta=10^{\circ}$. It however appears that for the best choice of homodyne angle, the classical-only estimator always outperforms the coherent-classical estimator.\vspace*{-2mm}

\begin{conjecture}\label{thm:central_result3}
Consider a coherent-classical estimation scheme defined by (\ref{eq:plant}), (\ref{eq:coh_controller}), (\ref{eq:coh_class_hd}) and (\ref{eq:coh_class_estimator}) with a cost $\tilde{J}_c$ defined in (\ref{eq:coh_cost}). Also, consider the corresponding purely-classical estimation scheme defined by (\ref{eq:plant}), (\ref{eq:class_hd}) and (\ref{eq:class_estimator}) with a cost $\bar{J}_c$ defined in (\ref{eq:class_cost}). Then, for the optimal choice of the homodyne angle $\theta_{opt}$,
\vspace*{-6mm}\small
\begin{equation}
\tilde{J}_c(\theta_{opt}) \geq \bar{J}_c(\theta_{opt}).
\end{equation}
\normalsize \vspace*{-6mm}
\end{conjecture}

\section{Coherent-Classical Estimation with Feedback}\label{sec:ccef}
\begin{figure}[!b]
\centering
\includegraphics[width=0.5\textwidth]{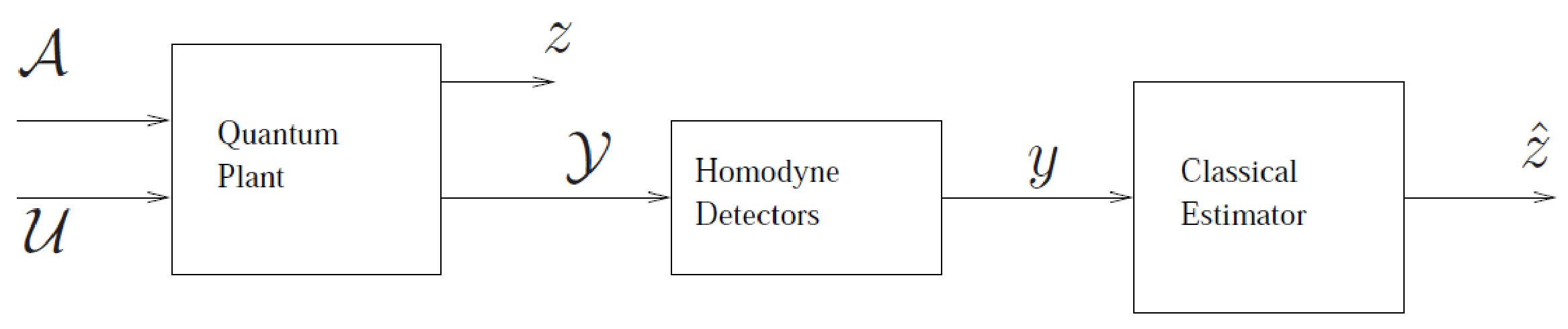}
\caption{\small Modified schematic of purely-classical estimation.}
\label{fig:fb_cls_scm}
\end{figure}

\vspace*{-3mm}
Here, we consider the case where there is quantum feedback from the coherent controller to the quantum plant \cite{IRP2}. For this purpose, the plant is assumed to have a control input $\mathcal{U}$ as in Fig. \ref{fig:fb_cls_scm}. Then, (\ref{eq:plant}) becomes
\vspace*{-7mm}\small
\begin{equation}\label{eq:fb_plant}
\begin{split}
\left[\begin{array}{c}
da\\
da^{\#}
\end{array}\right] &= F \left[\begin{array}{c}
a\\
a^{\#}
\end{array}\right] dt + \left[\begin{array}{cc}
G_1 & G_2
\end{array}\right] \left[\begin{array}{c}
d\mathcal{A}\\
d\mathcal{A}^{\#}\\
d\mathcal{U}\\
d\mathcal{U}^{\#}
\end{array}\right];\\
\left[\begin{array}{c}
d\mathcal{Y}\\
d\mathcal{Y}^{\#}
\end{array}\right] &= H \left[\begin{array}{c}
a\\
a^{\#}
\end{array}\right] dt + \left[\begin{array}{cc}
K & 0
\end{array}\right] \left[\begin{array}{c}
d\mathcal{A}\\
d\mathcal{A}^{\#}\\
d\mathcal{U}\\
d\mathcal{U}^{\#}
\end{array}\right];\\
z &= C\left[\begin{array}{c}
a\\
a^{\#} 
\end{array}\right].
\end{split}
\end{equation}
\normalsize \vspace*{-6mm}

The optimal purely-classical estimator is obtained from the solution of a Riccati equation of the form (\ref{eq:class_riccati}):
\vspace*{-6mm}\small
\begin{equation}\label{eq:fb_class_riccati}
\begin{split}
F\bar{P}_e &+ \bar{P}_eF^\dagger + G_1G_1^\dagger + G_2G_2^\dagger - (G_1 + \bar{P}_eH^\dagger)\\
&\times L^\dagger (LL^\dagger)^{-1} L(G_1 + \bar{P}_eH^\dagger)^\dagger =0,
\end{split}
\end{equation}
\normalsize \vspace*{-6mm}

where we have assumed $K=I$ as before. The estimation error cost is then given by (\ref{eq:class_cost}).\vspace*{-2mm}

\begin{figure}
\centering
\includegraphics[width=0.5\textwidth]{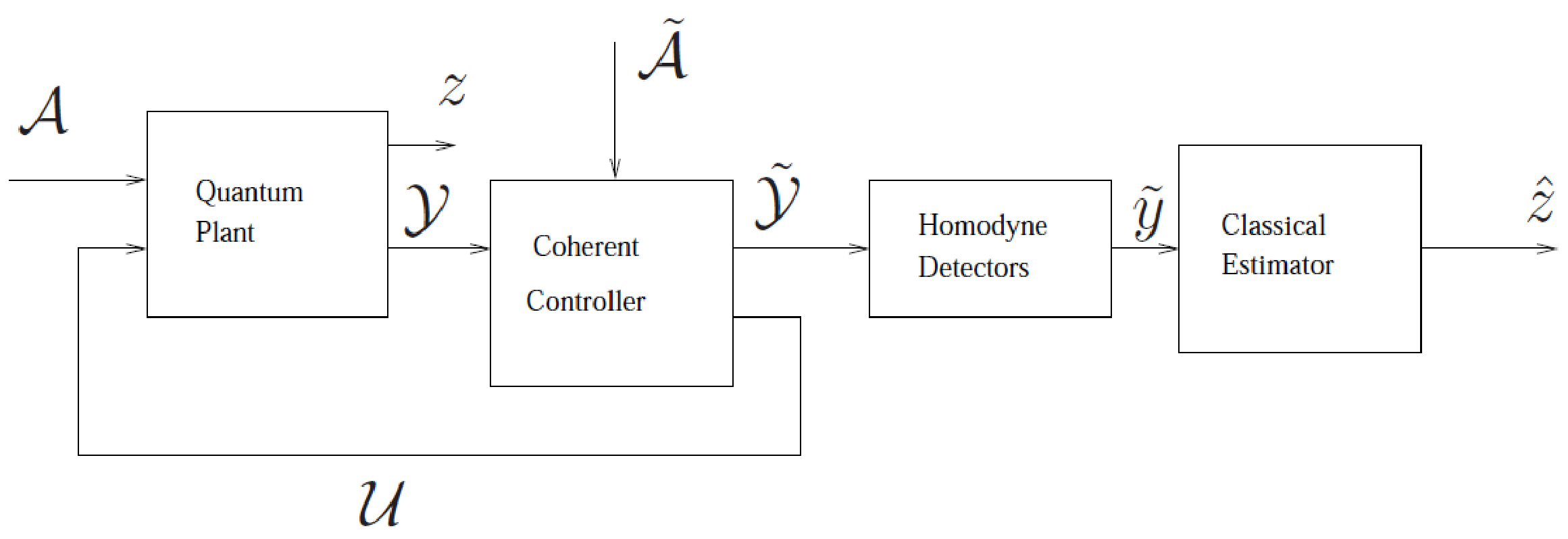}
\caption{\small Schematic diagram of coherent-classical estimation with coherent feedback.}
\label{fig:fb_coh_cls_scm}
\end{figure}

The coherent controller here would have an additional output that is fed back to the control input of the quantum plant, as depicted in Fig. \ref{fig:fb_coh_cls_scm}. The coherent controller in this case is defined as follows \cite{IRP2}:
\vspace*{-7mm}\small
\begin{equation}\label{eq:fb_coh_controller}
\begin{split}
\left[\begin{array}{c}
da_c\\
da_c^{\#}
\end{array}\right] &= F_c \left[\begin{array}{c}
a_c\\
a_c^{\#}
\end{array}\right] dt + \left[\begin{array}{cc}
G_{c1} & G_{c2}
\end{array}\right] \left[\begin{array}{c}
d\tilde{\mathcal{A}}\\
d\tilde{\mathcal{A}}^{\#}\\
d\mathcal{Y}\\
d\mathcal{Y}^{\#} 
\end{array}\right];\\
\left[\begin{array}{c}
d\tilde{\mathcal{Y}}\\
d\tilde{\mathcal{Y}}^{\#}\\
d\mathcal{U}\\
d\mathcal{U}^{\#}
\end{array}\right] &= \left[\begin{array}{c}
\tilde{H}_c\\
H_c
\end{array}\right] \left[\begin{array}{c}
a_c\\
a_c^{\#}
\end{array}\right] dt + \left[\begin{array}{cc}
\tilde{K}_{c1} & \tilde{K}_{c2}\\
K_{c1} & K_{c2}
\end{array}\right] \left[\begin{array}{c}
d\tilde{\mathcal{A}}\\
d\tilde{\mathcal{A}}^{\#}\\
d\mathcal{Y}\\
d\mathcal{Y}^{\#} 
\end{array}\right].
\end{split}
\end{equation}
\normalsize \vspace*{-6mm}

The plant (\ref{eq:fb_plant}) and the controller (\ref{eq:fb_coh_controller}) can be combined to yield an augmented system \cite{IRP2}:
\vspace*{-7mm}\small
\begin{equation}\label{eq:fb_coh_system}
\begin{split}
\left[\begin{array}{c}
da\\
da^{\#}\\
da_c\\
da_c^{\#}
\end{array}\right] &= \left[\begin{array}{cc}
F+G_2K_{c2}H & G_2H_c\\
G_{c2}H & F_c
\end{array}\right] \left[\begin{array}{c}
a\\
a^{\#}\\
a_c\\
a_c^{\#}
\end{array}\right] dt\\
&+ \left[\begin{array}{cc}
G_1+G_2K_{c2}K & G_2K_{c1}\\
G_{c2}K & G_{c1}
\end{array}\right] \left[\begin{array}{c}
d\mathcal{A}\\
d\mathcal{A}^{\#}\\
d\tilde{\mathcal{A}}\\
d\tilde{\mathcal{A}}^{\#}
\end{array}\right];\\
\left[\begin{array}{c}
d\tilde{\mathcal{Y}}\\
d\tilde{\mathcal{Y}}^{\#}
\end{array}\right] &= \left[\begin{array}{cc}
\tilde{K}_{c2}H & \tilde{H}_c
\end{array}\right] \left[\begin{array}{c}
a\\
a^{\#}\\
a_c\\
a_c^{\#}
\end{array}\right] dt + \left[\begin{array}{cc}
\tilde{K}_{c2}K & \tilde{K}_{c1}
\end{array}\right] \left[\begin{array}{c}
d\mathcal{A}\\
d\mathcal{A}^{\#}\\
d\tilde{\mathcal{A}}\\
d\tilde{\mathcal{A}}^{\#}
\end{array}\right].
\end{split}
\end{equation}
\normalsize \vspace*{-6mm}

The optimal coherent-classical estimator is then obtained from the solution $\tilde{P}_e$ (given by (\ref{eq:coh_err_cov})) to an algebraic Riccati equation of the form (\ref{eq:riccati}), where
\vspace*{-7mm}\small
\begin{equation}\label{eq:fb_coh_sys_matrices}
\begin{split}
F_a &= \left[\begin{array}{cc}
F+G_2K_{c2}H & G_2H_c\\
G_{c2}H & F_c
\end{array}\right],\\
G_a &= \left[\begin{array}{cc}
G_1+G_2K_{c2}K & G_2K_{c1}\\
G_{c2}K & G_{c1}
\end{array}\right],\\
H_a &= \left[\begin{array}{cc}
\tilde{K}_{c2}H & \tilde{H}_c
\end{array}\right], \qquad K_a = \left[\begin{array}{cc}
\tilde{K}_{c2}K & \tilde{K}_{c1}
\end{array}\right],
\end{split}
\end{equation}
\normalsize \vspace*{-6mm}

and $\tilde{L}_1$, $\tilde{L}_2$ and $L$ as in (\ref{eq:coh_sys_matrices}). Here, for the coherent controller to be physically realizable, we would have:
\vspace*{-6mm}\small \[ \left[\begin{array}{cc}
\tilde{K}_{c1} & \tilde{K}_{c2}\\
K_{c1} & K_{c2}
\end{array}\right] = I, \]
\normalsize \vspace*{-5mm}

which implies $\tilde{K}_{c1}=K_{c2}=I$ and $K_{c1}=\tilde{K}_{c2}=0$. The estimation error is then given by the cost (\ref{eq:coh_cost}).\vspace*{-2mm}

\begin{remark}
Note that the combined plant-controller system being measured here is again a fully quantum system. The coherent controller not only preserves the quantum coherence of the quantum plant output, but also allows for coherent feedback control of the quantum plant by means of a suitable choice of the controller parameters that minimizes (\ref{eq:coh_cost}). This further assists in improving the precision of the classical estimate of a plant variable, when compared to the cases of purely-classical estimation and coherent-classical estimation without coherent feedback.
\end{remark}

\begin{theorem}\label{thm:central_result4}
Consider a coherent-classical estimation scheme defined by (\ref{eq:fb_plant}), (\ref{eq:fb_coh_controller}), (\ref{eq:coh_class_hd}) and (\ref{eq:coh_class_estimator}), such that both the plant and the controller are physically realizable annihilation operator only systems, with the cost $\tilde{J}_c$ as in (\ref{eq:coh_cost}). Also, consider the corresponding purely-classical estimation scheme defined by (\ref{eq:fb_plant}), (\ref{eq:class_hd}) and (\ref{eq:class_estimator}), such that the plant is a physically realizable annihilation operator only system, with the cost $\bar{J}_c$ as in (\ref{eq:class_cost}). Then,
\vspace*{-7mm}\small
\begin{equation}
\tilde{J}_c = \bar{J}_c.
\end{equation}
\normalsize \vspace*{-5mm}
\end{theorem}
\vspace*{-2mm}
\begin{proof}
The plant (\ref{eq:fb_plant}) may be augmented to account for an unused output $\bar{\mathcal{Y}}$ to recast the QSDE's in the desired form, that lends itself appropriately to the physical realizability treatment, as follows:
\vspace*{-1.1cm}\small
\begin{equation}\label{eq:fb_mod_plant}
\begin{split}
\left[\begin{array}{c}
da\\
da^{\#}
\end{array}\right] &= F \left[\begin{array}{c}
a\\
a^{\#}
\end{array}\right] dt + \left[\begin{array}{cc}
G_1 & G_2
\end{array}\right] \left[\begin{array}{c}
d\mathcal{A}\\
d\mathcal{A}^{\#}\\
d\mathcal{U}\\
d\mathcal{U}^{\#}
\end{array}\right];\\
\left[\begin{array}{c}
d\mathcal{Y}\\
d\mathcal{Y}^{\#}\\
d\bar{\mathcal{Y}}\\
d\bar{\mathcal{Y}}^{\#}
\end{array}\right] &= \left[\begin{array}{c}
H\\
\bar{H}
\end{array}\right] \left[\begin{array}{c}
a\\
a^{\#}
\end{array}\right] dt + \left[\begin{array}{cc}
K & 0\\
0 & \bar{K}
\end{array}\right] \left[\begin{array}{c}
d\mathcal{A}\\
d\mathcal{A}^{\#}\\
d\mathcal{U}\\
d\mathcal{U}^{\#}
\end{array}\right];\\
z &= C\left[\begin{array}{c}
a\\
a^{\#} 
\end{array}\right].
\end{split}
\end{equation}
\normalsize \vspace*{-6mm}

Here, $\bar{K}=I$ for the plant to be physically realizable. Additionally, using the same arguments as in the proof for Theorem \ref{thm:central_result1}, we must have:
\vspace*{-7mm}\small
\begin{equation}\label{eq:fb_plant_phys_rlz}
\begin{split}
F\Theta +\Theta F^\dagger +G_1G_1^\dagger +G_2G_2^\dagger &=0,\\
G_1 &=-\Theta H^\dagger ,\\
G_2 &=-\Theta \bar{H}^\dagger ,
\end{split}
\end{equation}
\normalsize \vspace*{-6mm}

for the annihilation operator only plant to be physically realizable with the commutation matrix $\Theta > 0$.\vspace*{-2mm}

Similarly, if the coherent controller (\ref{eq:fb_coh_controller}) is an annihilation operator only system, we must have the following for it to be physically realizable:
\vspace*{-7mm}\small
\begin{equation}\label{eq:fb_ctrl_phys_rlz}
\begin{split}
F_c\Theta_c +\Theta_c F_c^\dagger +G_{c1}G_{c1}^\dagger +G_{c2}G_{c2}^\dagger &=0,\\
G_{c1} &=-\Theta_c \tilde{H}_c^\dagger ,\\
G_{c2} &=-\Theta_c H_c^\dagger ,
\end{split}
\end{equation}
\normalsize \vspace*{-6mm}

where $\Theta_c > 0$ is the controller's commutation matrix.\vspace*{-2mm}

Clearly, $\bar{P}_e = \Theta$ satisfies the Riccati equation (\ref{eq:fb_class_riccati}), owing to (\ref{eq:fb_plant_phys_rlz}), for the purely-classical estimation case. Moreover, it follows from (\ref{eq:fb_plant_phys_rlz}) and (\ref{eq:fb_ctrl_phys_rlz}), that $\tilde{P}_e = \left[\begin{array}{cc}
\Theta & 0\\
0 & \Theta_c
\end{array}\right]$ satisfies 
(\ref{eq:riccati}), (\ref{eq:fb_coh_sys_matrices}). Thus, we get $\bar{J}_c = \tilde{J}_c = C\Theta C^\dagger$.
\hfill \qed
\end{proof}

\begin{remark}
Theorem \ref{thm:central_result4} implies that coherent-classical estimation with coherent feedback, where both the plant and the controller are physically realizable annihilation operator quantum systems, performs identical to, and no better than, purely-classical estimation of the plant. Note that in addition to $P_2=0$, we need to have both $P_1=\Theta$ and $P_3=\Theta_c$ for the coherent-classical scheme to be equivalent to the classical-only scheme.
\end{remark}

Now, we present examples involving dynamic squeezers for the case of coherent-classical estimation with feedback. First, we give one to illustrate Theorem \ref{thm:central_result4}. Here, the quantum plant (\ref{eq:sqz_plant}) takes the form:
\vspace*{-8mm}\small
\begin{equation}\label{eq:fb_sqz_plant}
\begin{split}
\left[\begin{array}{c}
da\\
da^{*}
\end{array}\right] &= \left[\begin{array}{cc}
-\frac{\gamma}{2} & -\chi\\
-\chi^{*} & -\frac{\gamma}{2}
\end{array}\right] \left[\begin{array}{c}
a\\
a^{*}
\end{array}\right] dt\\
&- \sqrt{\kappa_1} \left[\begin{array}{c}
d\mathcal{A}\\
d\mathcal{A}^{*}
\end{array}\right] - \sqrt{\kappa_2} \left[\begin{array}{c}
d\mathcal{U}\\
d\mathcal{U}^{*}
\end{array}\right];\\
\left[\begin{array}{c}
d\mathcal{Y}\\
d\mathcal{Y}^{*}
\end{array}\right] &= \sqrt{\kappa} \left[\begin{array}{c}
a\\
a^{*}
\end{array}\right] dt + \left[\begin{array}{c}
d\mathcal{A}\\
d\mathcal{A}^{*}
\end{array}\right];\\
z &= \left[\begin{array}{cc}
\frac{1}{\sqrt{2}} & -\frac{1}{\sqrt{2}}
\end{array}\right] \left[\begin{array}{c}
a\\
a^{*}
\end{array}\right].
\end{split}
\end{equation}
\normalsize \vspace*{-6mm}

Here, we choose $\gamma = 4$, $\kappa_1 = \kappa_2 = 2$ and $\chi = 0$. Note that this system is physically realizable, since $\gamma = \kappa_1+\kappa_2$, and is annihilation operator only, since $\chi = 0$. We then calculate the optimal classical-only state estimator and the error $\bar{J}_c$ in (\ref{eq:class_cost}) for this system using the standard Kalman filter equations corresponding to homodyne detector angles varying from $\theta = 0^{\circ}$ to $\theta = 180^{\circ}$.\vspace*{-2mm}

The coherent controller (\ref{eq:sqz_ctrlr}) in this case takes the form:
\vspace*{-7mm}\small
\begin{equation}\label{eq:fb_sqz_ctrlr}
\begin{split}
\left[\begin{array}{c}
da\\
da^{*}
\end{array}\right] &= \left[\begin{array}{cc}
-\frac{\gamma}{2} & -\chi\\
-\chi^{*} & -\frac{\gamma}{2}
\end{array}\right] \left[\begin{array}{c}
a\\
a^{*}
\end{array}\right] dt\\
&- \sqrt{\kappa_1} \left[\begin{array}{c}
d\tilde{\mathcal{A}}\\
d\tilde{\mathcal{A}}^{*}
\end{array}\right] - \sqrt{\kappa_2} \left[\begin{array}{c}
d\mathcal{Y}\\
d\mathcal{Y}^{*}
\end{array}\right];\\
\left[\begin{array}{c}
d\tilde{\mathcal{Y}}\\
d\tilde{\mathcal{Y}}^{*}
\end{array}\right] &= \sqrt{\kappa_1} \left[\begin{array}{c}
a\\
a^{*}
\end{array}\right] dt + \left[\begin{array}{c}
d\tilde{\mathcal{A}}\\
d\tilde{\mathcal{A}}^{*}
\end{array}\right];\\
\left[\begin{array}{c}
d\mathcal{U}\\
d\mathcal{U}^{*}
\end{array}\right] &= \sqrt{\kappa_2} \left[\begin{array}{c}
a\\
a^{*}
\end{array}\right] dt + \left[\begin{array}{c}
d\mathcal{Y}\\
d\mathcal{Y}^{*}
\end{array}\right].
\end{split}
\end{equation}
\normalsize \vspace*{-6mm}

\begin{figure}[!t]
\centering
\includegraphics[width=0.5\textwidth]{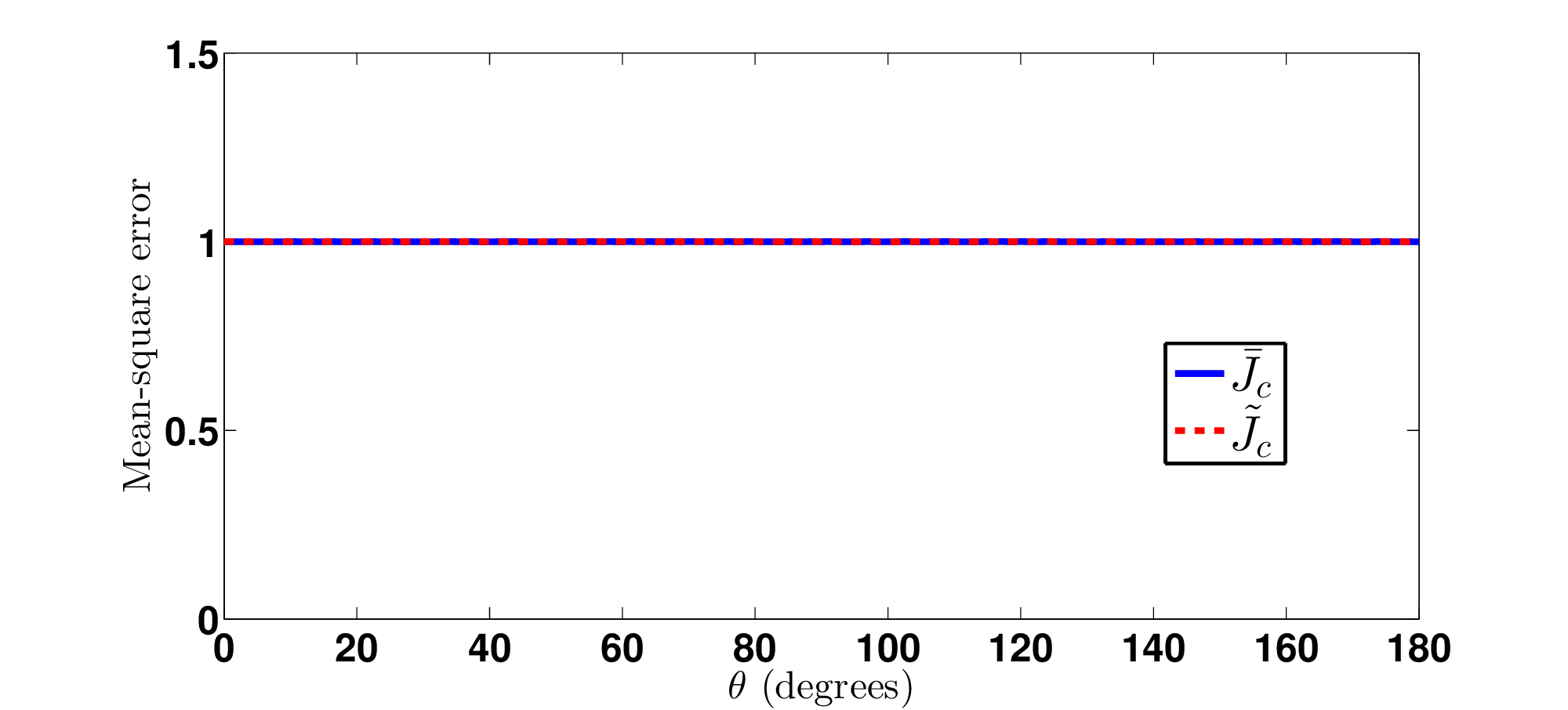}
\caption{\small Feedback: Estimation error vs. homodyne angle $\theta$ in the case of annihilation operator only plant and controller.}
\label{fig:cls_vs_coh_cls_4}
\end{figure}

\begin{figure}[!b]
\centering
\includegraphics[width=0.5\textwidth]{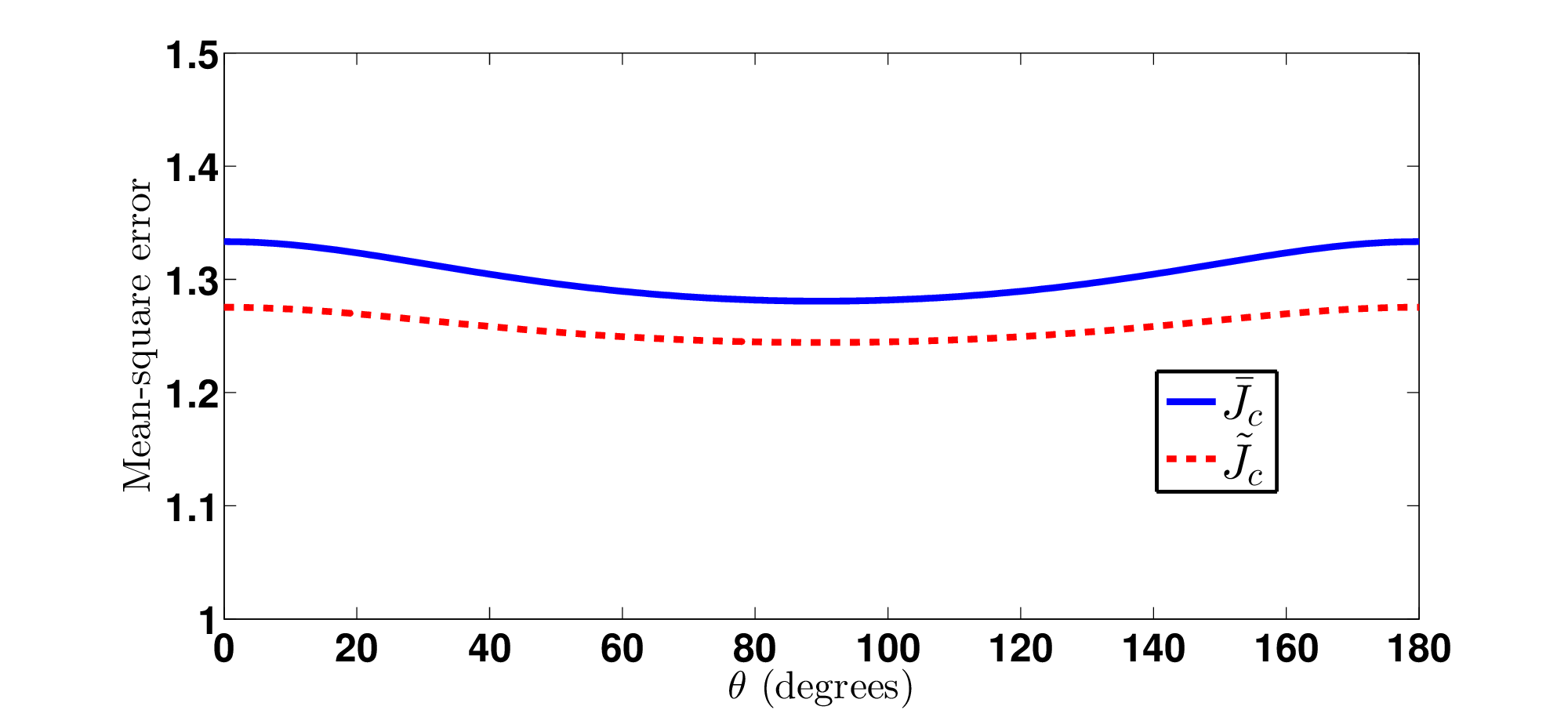}
\caption{\small Feedback: Estimation error vs. homodyne angle $\theta$ in the case of an annihilation operator only controller.}
\label{fig:cls_vs_coh_cls_5}
\end{figure}

Here, we choose $\gamma = 16$, $\kappa_1 = \kappa_2 = 8$ and $\chi = 0$, so that it is a physically realizable annihilation operator only system. Then, the classical estimator for this case is calculated according to (\ref{eq:fb_coh_system}), (\ref{eq:fb_coh_sys_matrices}), (\ref{eq:riccati}), (\ref{eq:coh_sys_est}) for the homodyne detector angle varying from $\theta=0^{\circ}$ to $\theta=180^{\circ}$. The resulting value of the cost $\tilde{J}_c$ in (\ref{eq:coh_cost}) alongwith the cost for the purely-classical estimator case is shown in Fig. \ref{fig:cls_vs_coh_cls_4}. Clearly both the classical-only and coherent-classical estimators have the same estimation error cost for all homodyne angles. This illustrates Theorem \ref{thm:central_result4}.\vspace*{-2mm}

We now show in examples that when either the plant or the controller is not an annihilation operator quantum system, the coherent-classical estimator with coherent feedback can provide improvement over the purely-classical estimator. We first consider an example where the controller is a physically realizable annihilation operator only system. But, the plant is physically realizable with $\chi \neq 0$. In (\ref{eq:fb_sqz_plant}), we choose $\gamma = 4$, $\kappa_1 = \kappa_2 = 2$, $\chi = 0.5$, and in (\ref{eq:fb_sqz_ctrlr}), $\gamma = 16$, $\kappa_1 = \kappa_2 = 8$, $\chi = 0$. Fig. \ref{fig:cls_vs_coh_cls_5} then shows that the coherent-classical error is less than the purely-classical error for all homodyne angles.\vspace*{-2mm}

\begin{figure}[!t]
\centering
\includegraphics[width=0.5\textwidth]{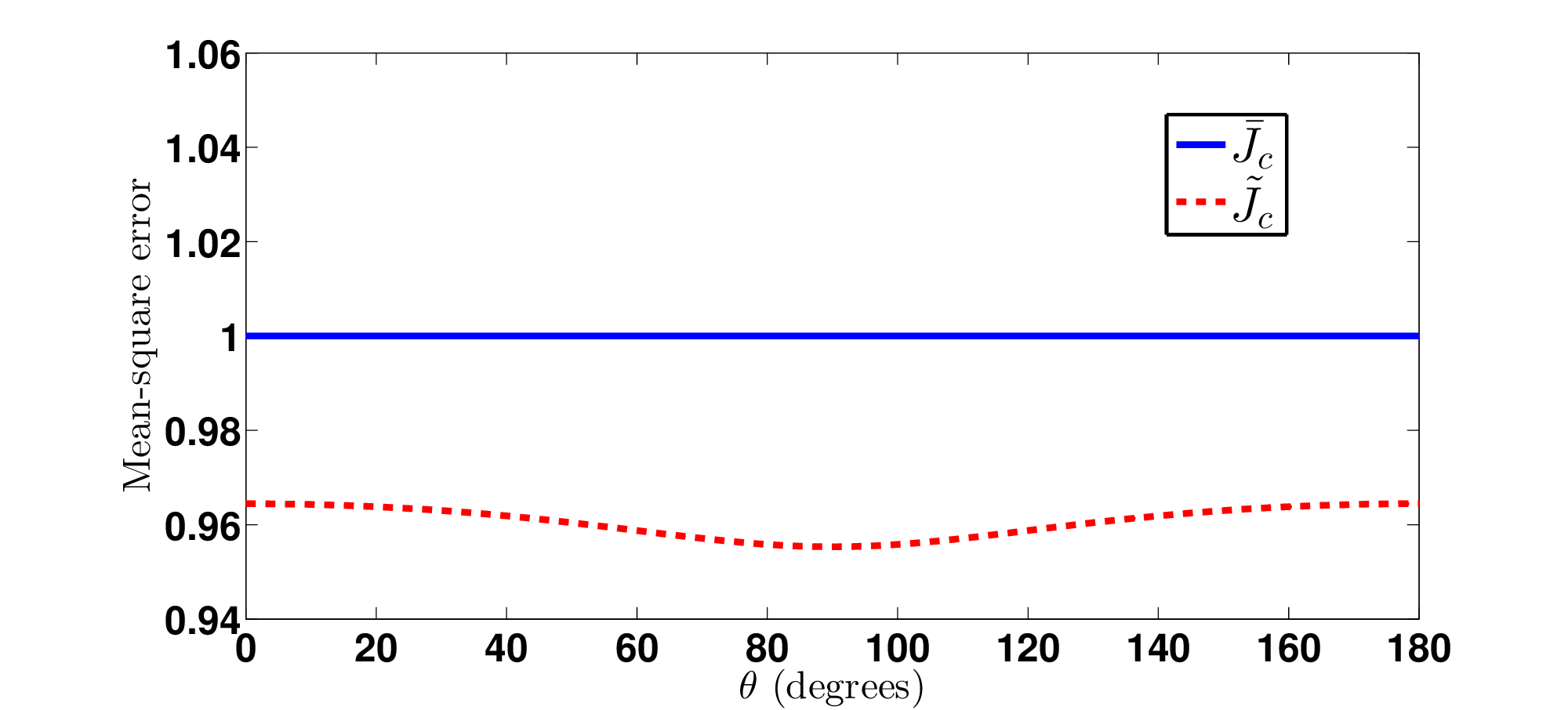}
\caption{\small Feedback: Estimation error vs. homodyne angle $\theta$ in the case of an annihilation operator only plant.}
\label{fig:cls_vs_coh_cls_6}
\end{figure}

Then, we consider the case where the plant is an annihilation operator only system, but the coherent controller is not. In (\ref{eq:fb_sqz_plant}), we choose $\gamma = 4$, $\kappa_1 = \kappa_2 = 2$, $\chi = 0$, and in (\ref{eq:fb_sqz_ctrlr}), we take $\gamma = 16$, $\kappa_1 = \kappa_2 = 8$, $\chi = -0.5$. Fig. \ref{fig:cls_vs_coh_cls_6} then shows that the coherent-classical error is again less than the purely-classical error for all homodyne angles.\vspace*{-2mm}

\begin{figure}[!t]
\centering
\includegraphics[width=0.5\textwidth]{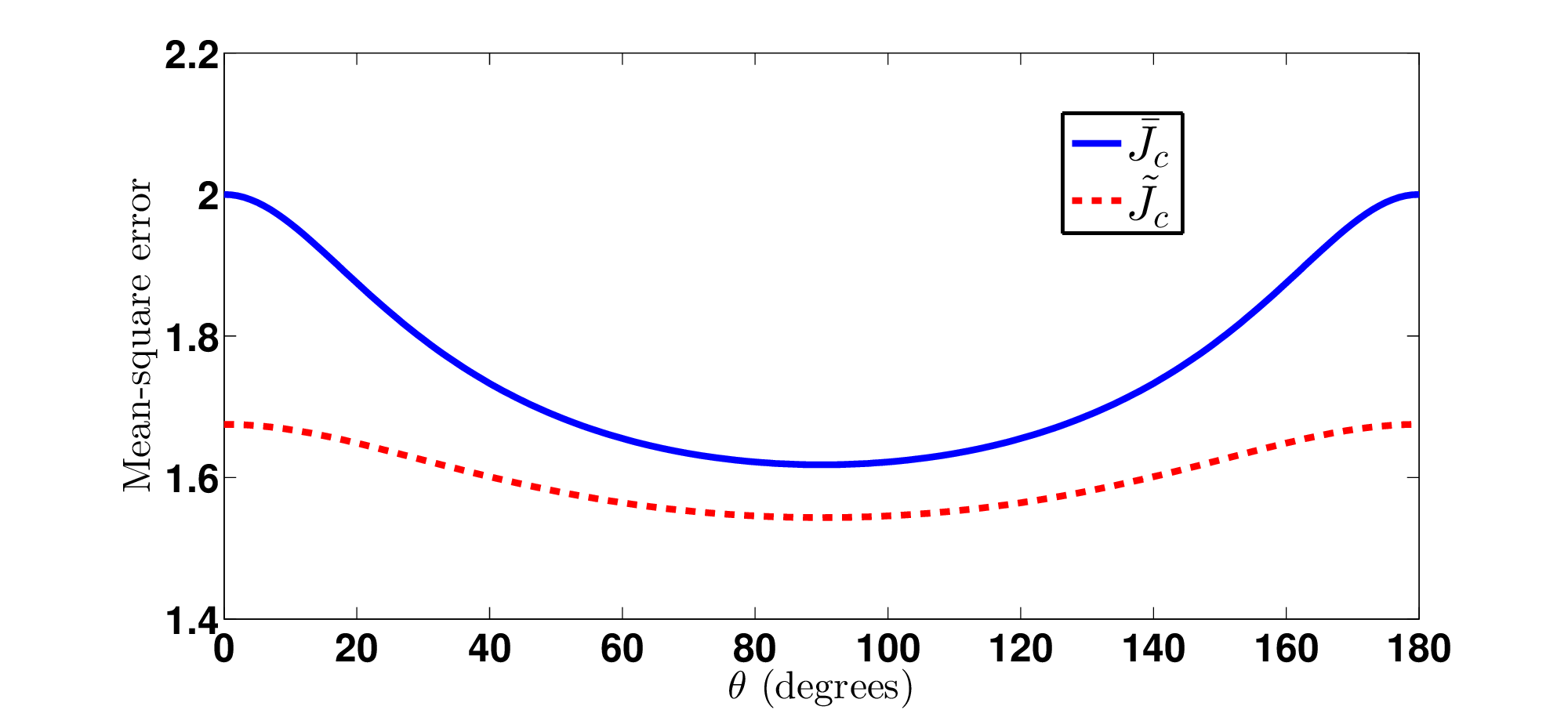}
\caption{\small Feedback: Estimation error vs. homodyne angle $\theta$ in the case of a squeezer plant and a squeezer controller.}
\label{fig:cls_vs_coh_cls_7}
\end{figure}

\begin{figure}[!b]
\centering
\includegraphics[width=0.5\textwidth]{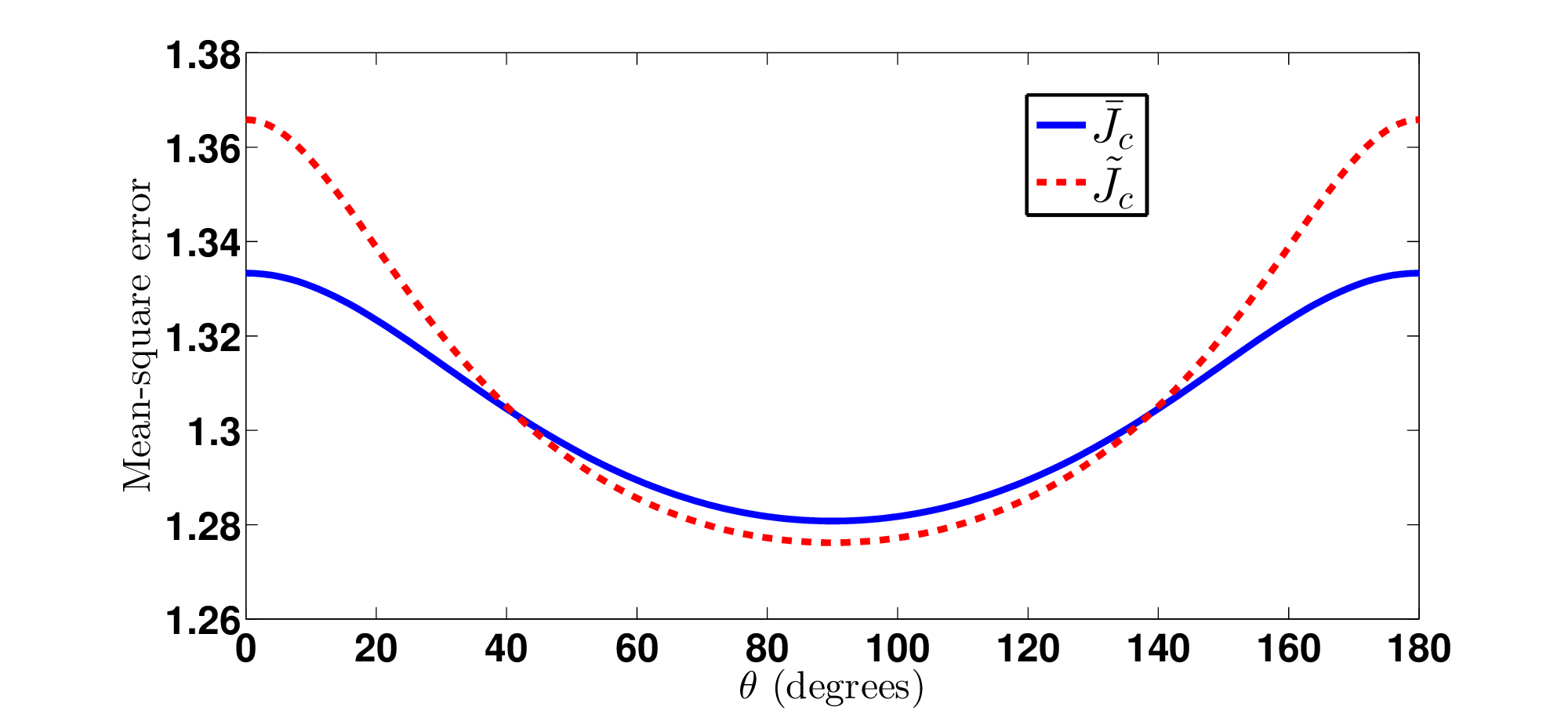}
\caption{\small Feedback: Estimation error vs. homodyne angle $\theta$ in the case of a squeezer plant and a squeezer controller, where it is possible to get better coherent-classical estimates than purely-classical estimates only for certain homodyne angles.}
\label{fig:cls_vs_coh_cls_8}
\end{figure}

We also show the case where both the plant and the controller have $\chi \neq 0$. With $\gamma = 4$, $\kappa_1 = \kappa_2 = 2$, $\chi = 1$ in (\ref{eq:fb_sqz_plant}), and $\gamma = 16$, $\kappa_1 = \kappa_2 = 8$, $\chi = -0.5$ in (\ref{eq:fb_sqz_ctrlr}), Fig. \ref{fig:cls_vs_coh_cls_7} shows that coherent-classical error is less than purely-classical error for all homodyne angles. However, with both plant and controller having $\chi \neq 0$, coherent-classical estimates can be better than purely-classical estimates for only certain homodyne angles, as in Fig. \ref{fig:cls_vs_coh_cls_8}, where we used $\gamma = 4$, $\kappa_1 = \kappa_2 = 2$, $\chi = 0.5$ in (\ref{eq:fb_sqz_plant}) and $\gamma = 16$, $\kappa_1 = \kappa_2 = 8$, $\chi = 0.5$ in (\ref{eq:fb_sqz_ctrlr}).\vspace*{-2mm}

We observe that if there is any improvement with the coherent-classical estimation (with feedback) over purely-classical estimation, the former is always superior to the latter for the best choice of the homodyne angle. This we propose as a conjecture here. This is just the opposite of Conjecture \ref{thm:central_result3} for the no feedback case.\vspace*{-2mm}

\begin{conjecture}\label{thm:central_result5}
Consider a coherent-classical estimation scheme defined by (\ref{eq:fb_plant}), (\ref{eq:fb_coh_controller}), (\ref{eq:coh_class_hd}) and (\ref{eq:coh_class_estimator}) with a cost $\tilde{J}_c$ defined in (\ref{eq:coh_cost}). Also, consider the corresponding purely-classical estimation scheme defined by (\ref{eq:fb_plant}), (\ref{eq:class_hd}) and (\ref{eq:class_estimator}) with a cost $\bar{J}_c$ defined in (\ref{eq:class_cost}). Then, if there exists a homodyne angle $\theta_i$ for which  $\tilde{J}_c(\theta_i) \leq \bar{J}_c(\theta_i)$, for the best choice $\theta_{opt}$ of homodyne angle,
\vspace*{-6mm}\small
\begin{equation}
\tilde{J}_c(\theta_{opt}) \leq \bar{J}_c(\theta_{opt}).
\end{equation}
\normalsize \vspace*{-6mm}
\end{conjecture}

\vspace*{-2mm}All of the above results in this section suggest that either or both of the plant and the controller need to have non-zero squeezing to produce better estimates than in classical-only case. This is because the coherent information (in the spirit of Ref. \cite{SN}) at the output of the combined plant-controller quantum system here can be no more than at the output of the quantum plant, when both the plant and the controller are passive systems.

\vspace*{-3mm}
\section{Conclusion}\label{sec:conc}\vspace*{-3mm}
In this paper, we studied two flavours of coherent-classical estimator, one with coherent feedback and the other without, for a class of linear quantum systems. We did a comparison study of these with the corresponding purely-classical estimators. Indeed, the class of linear quantum systems considered here can be rewritten in terms of linear classical stochastic systems and the results explained in the classical world. However, any classical model for the plant or the controller is inherently physical, whereas the corresponding models depicting a quantum plant or controller need to satisfy the physical realizability constraints to be actual physical systems. Moreover, our results imply physically that the combined plant-controller quantum system under measurement should have increased coherent information at its output compared with the output of the plant alone being measured, to be able to produce more accurate classical estimates of a plant variable. Intuitively, these results should also hold for non-linear quantum systems.

\vspace*{-3mm}
\begin{ack}\vspace*{-3mm}
This work was supported by the Australian Research Council (ARC) under grants CE110001027 (EHH) and FL110100020 (IRP), and by the US Air Force Office of Scientific Research (AFOSR) under agreement number FA2386-16-1-4065 (IRP). SR was funded by the Singapore National Research Foundation Grant No.~NRF-NRFF2011-07 and the Singapore Ministry of Education Academic Research Fund Tier 1 Project R-263-000-C06-112, and is currently funded by the UK National Quantum Technologies Programme (EP/M01326X/1, EP/M013243/1). Moreover, SR thanks Mohamed Mabrok for useful discussions related to this work.
\end{ack}
\vspace*{-3mm}
\bibliographystyle{ifacconf}
\bibliography{cohclsbib}

\end{document}